\documentclass[12 pt, a4paper, twoside, openright]{article}

\usepackage{amsmath, amssymb, theorem, latexsym}

\usepackage[english]{babel}
\newtheorem{theorem}{Theorem}[section]
\newtheorem{lemma}[theorem]{Lemma}

\newtheorem{proposition}[theorem]{Proposition}
\newtheorem{definition}[theorem]{Definition}

\newtheorem{corollary}[theorem]{Corollary}

\newtheorem{remark}[theorem]{Remark}

\allowdisplaybreaks
\renewcommand{\thefootnote}{\fnsymbol{footnote}}

\numberwithin{equation}{section}

\title {\textsf {Positive solutions of Schr\"odinger equations and fine regularity of boundary points} }

\author { Alano Ancona{$^{}$}\\ {\small D\'{e}partement de Math\'{e}matiques, Universit\'{e} Paris-Sud 11}
\\ {\small Orsay 91\hspace{0.5pt}405 France}} 
\date{\empty}
\date{{\small   March 6, 2012 Version 3\footnote{Mainly an annoying typo in Theorem 1.1 is corrected}
}}

\begin{document}

\maketitle

{\renewcommand\thefootnote{}
\footnote{1991 Mathematics Subject Classification 31C35, 31C25, 35J15, 35C15
}}

\vspace{10truemm}

\renewcommand{\abstractname}

\noindent{\bf Abstract.} {\footnotesize Given a Lipschitz domain $\Omega $ in ${ \mathbb R} ^N $  and a nonnegative potential $V$ in $\Omega $ such that $V(x)\, d(x,\partial \Omega )^2$ is bounded  we study the {\sl fine regularity} of  boundary points  with respect to the  Schr\"odinger  operator $L_V:= \Delta  -V$ in $\Omega $.  Using potential theoretic methods, several  conditions are shown to be equivalent to the fine regularity  of $z \in   \partial \Omega $. The main result is a  simple ({\em explicit} if $\Omega $ is smooth) necessary and sufficient condition involving  the size of $V$ for $z \in   \partial \Omega $ to be finely regular. An intermediate result consists in a majorization of $\int_A  \vert  { \frac {  u} {d(.,\partial \Omega ) }} \vert  ^2\, dx$ for $u$ positive harmonic in $\Omega $ and $A \subset   \Omega $. Conditions for almost everywhere regularity in a subset $A $ of $   \partial \Omega $ are also given as well as an extension of the main results to a notion of fine ${ \mathcal L}_1 \vert  { \mathcal L}_0$-regularity, if ${ \mathcal L}_j={ \mathcal L}-V_j$, $V_0,\, V_1$ being two potentials, with $V_0 \leq V_1$ and ${ \mathcal L}$ a second order elliptic operator.}

\noindent

\section{Introduction and main results.}This paper stems from  a question raised by Moshe Marcus and Laurent V\'{e}ron several years ago and its answer   \cite{anc4} (see the appendix in \cite{VYA} and the comments after Theorem \ref{charac} below). Another independent -- and recent -- question of Moshe Marcus has also motivated the approach followed in Section 3.  See Proposition \ref{MM}.

Both questions dealt with positive solutions  of a Schr\" odinger equation  $ \Delta  u-Vu=0$ in a  Lipschitz domain $\Omega $ with $V$ in a natural class of nonnegative potentials in $\Omega $. The first was about a necessary explicit condition for a boundary point of  $\Omega $ to be ``finely regular" with respect to $\Delta  -V$ and is related to the works \cite{MV1} \cite{MV5}.  The definition of ``fine regularity" \cite{VYA} recalled below goes back to a notion  introduced by E.\ B.\ Dynkin  to study  the  boundary values (or traces on $\partial \Omega $) of  positive solutions  of  nonlinear equation such as $ \Delta  u= u^ \alpha $, $ \alpha >1$ -- in which case, given a solution $u$,  Dynkin's definition corresponds here to $V= \vert  u \vert  ^{ \alpha -1}$ -- (see \cite{dyn0}, \cite{dynkuz} and the books \cite{dyn1}, \cite{dyn2}).

To state our main results we fix  some notations and recall basic definitions and facts  which will be used all along this paper.  This will be quite similar to  the appendix of \cite{VYA}. Other notations (that also appear in \cite{VYA}) are fixed later  in Section  2.2.

 Let  $\Omega $ be a bounded Lipschitz domain in ${ \mathbb R} ^N$, $N\geq 2$. For $x \in   \Omega $,   $ \delta_{\Omega } (x)$ is
the distance from $x$ to ${  \mathbb  R}^n\setminus \Omega$ in ${ \mathbb R} ^N$ and for $a \geq 0$, ${ \mathcal V}(\Omega ,a)$ (or ${ \mathcal V}_a(\Omega )$)  is the set of all nonnegative Borel measurable functions $V:\Omega \to  { \mathbb R} $ such that $V(x) \leq a/(\delta _{\Omega }(x))^2$ in $\Omega $. Let ${ \mathcal V}(\Omega )= \bigcup    _{a>0} { \mathcal V}(\Omega ,a)$. We also fix a  reference point $x_0$ in $\Omega $. 

For  $V \in   { \mathcal V}(\Omega)$, let $L_V:= \Delta  -V$  denote  the  corresponding Schr\" odinger  operator where $ \Delta  $ is the classical Laplacian. An $L_V$-harmonic function in $\Omega $ is the continuous representative of a weak solution in $\Omega $. It  belongs to $W^{2,p}_{\rm loc} (\Omega )$ for all $p< \infty $.

\noindent 

By the results in \cite{anc3} (see more details in Section 2)  it is  known that given $y \in   \partial \Omega $ and $V \in   { \mathcal V}(\Omega )$, there exists a unique positive  $L_V$-harmonic function $K_y^V$ in $\Omega $ that vanishes on $\partial \Omega \setminus  \{ y \}$ and satisfies $K_y^V(x_0)=1$. Moreover $y\mapsto K_y^V$ is continuous (for uniform convergence on compact subsets of $\Omega $). For $V=0$, we recover the well-known results of Hunt and Wheeden \cite{HW1},\cite{HW2} and we denote $K_y:=K_y^0$. 
The function $K^V_y$ is the Martin function with respect to $L_V$  in $\Omega $  with pole at $y$  normalized at $x_0$.  For general facts about Martin's boundary see e.g.\ \cite{martin,naim, doob, ancbook}.

The function $K_y$ is  superharmonic in $\Omega $ with respect to $L_V$ and we may hence consider its  greatest $L_V$-harmonic minorant $ \tilde K_y^V$ in $\Omega $. 
Clearly $ \tilde K_y^V=c(y)\, K_y^V$ with $c(y)= \tilde K_y^V(x_0)$, $0 \leq c(y) \leq 1$. Moreover $y\mapsto c(y)$ is upper semicontinuous.

\noindent

The set of ``{\em finely}" regular boundary points (ref.\ \cite{dyn1} Chap.\ 11, \cite {dynkuz}, \cite{VYA})   with respect to $L_V $ in $\Omega $ is 
  \begin{align}{ \mathcal R}\mbox{\footnotesize{eg}}_{_ V}(\Omega )&= \{ y \in   \partial \Omega\,;\,  \tilde K_y^V>0 \}= \{ y \in   \partial \Omega \,;\, c(y)>0\,  \}   \end{align}
\noindent and the set of finely irregular boundary points is ${ \mathcal S}\mbox{\footnotesize {ing}}_{_ V}(\Omega ):= \partial \Omega\setminus { \mathcal R}\mbox{\footnotesize{eg}}_{_ V}(\Omega )$. We  also say that $y \in   \partial \Omega $ is a $V$-finely regular boundary point if $y \in  { \mathcal R}\mbox{\footnotesize{eg}}_V(\Omega )$. Clearly   ${ \mathcal R}\mbox{\footnotesize{eg}}_{_ V}(\Omega )$ is a $K_  \sigma   $ subset of $\partial \Omega $. The probabilistic interpretation, which will not be used here, is that a   point $y \in   \partial \Omega$ is $V $-finely regular iff $\int _0^ \tau \, V( \xi  _s)\, ds <+ \infty $ a.s. for a Brownian motion  $ \{  \xi  _s \}_{0 \leq s< \tau }$ starting  at $x_0$ and conditioned to exit from $\Omega$ at  $y$. 

Denote $G^V(x,y)$ the Green's function with respect to $L_V$ in $\Omega $ and let $G:=G^0$. Thus for $y \in   \Omega $,   $G^V_y: x\mapsto G^V(x,y)$ is the smallest positive function in $\Omega $ which is $L_V$ harmonic in $\Omega \setminus  \{ y \}$ and such that   $G^V_y(y)=+ \infty $ and  $L_V(G_y)=-  \delta  _y$ in the weak sense. Here and below  $  \delta _y$ denotes  Dirac measure at $y$. 

Our main result about fine regularity is  as follows (see generalizations in Section 7). 

\begin{theorem} \label{th1}Let $V \in   { \mathcal V}(\Omega ,a)$. Given $y \in   \partial \Omega $, the following are equivalent:
\begin{enumerate}
\item The point $y$ is finely regular with respect to the potential $V$ in $\Omega $.

\item The integral $ \int_\Omega  G(x_0,z)\, V(z)\, K_y^V(z)\, dz$ is finite.

\item The integral $ \int_\Omega  G(x_0,z)\, V(z)\, K_y(z)\, dz$ is finite.
\end{enumerate}
\end{theorem}

For $\Omega $ smooth -- and using the well known fact that  away from $x_0$ the first eigenfunction of the Laplacian in $\Omega $ is equivalent in size to $G(x_0,.)$ -- the implication $(i)\Rightarrow(iii)$ appears also in \cite{VYA} and  the equivalence $(i)\Leftrightarrow (iii)$ is stated as an open question there. This implication is also clear using  the probabilistic point of view and the main point in Theorem  \ref{th1} is that $(iii)\Rightarrow(i)$.
\vspace{2truemm}

It follows from Theorem \ref{th1} that  if $\Omega $ is $C^{1, \alpha }$-smooth for some $ \alpha  \in   (0,1]$  then $y \in   \partial \Omega $ is $V$-finely regular if and only if 
$$\int _\Omega {\frac  {  \delta _\Omega (z)^2} { \vert  z-y \vert  ^N}}\, V(z)\, dz  < \infty .$$

\vspace{3truemm}

To establish Theorem \ref{th1}, we prove -- using methods in \cite{ancbook0} -- a general fact from potential theory in Lipschitz domains interesting in its own right. See Corollaries \ref{maincorogene} and \ref{maincorogene2}.

\begin{corollary} \label{cormain}Let $\mu  $ be a finite positive measure on $\partial \Omega $. Then $\mu  $-almost every point $y \in   \partial \Omega $ is $V$-finely regular iff there is a positive measure $\nu =f\cdot \mu  $, with $0<f \leq 1$, such that $\int_\Omega  G(x,x_0)\,V(x)\,K_\nu (x)\, dx< \infty $. If  $\Omega $ is $C^{1, \alpha }$-smooth, $0< \alpha  \leq 1$, this amounts to
 \begin{align} \iint_{\Omega \times\partial \Omega } {\frac  {  \delta _\Omega (x)^2} { \vert  x-y \vert  ^N}}\, V(x)\,\, dx\,d \nu (y)\, < +\infty . 
  & \nonumber \end{align}
\end{corollary}

In the next statement,  $H_{N-1}$ denotes the $N-1$-dimensional Hausdorff measure in ${ \mathbb  R}^N$.

\begin{theorem} \label{th2} Let $A$ be a Borel subset of $\partial \Omega $. Then, $H_{N-1}$-almost every point $y \in  A$ is finely regular (with respect to $\Omega $ and $V$) if and only if there exists, for $H_{N-1}$-almost every point $y \in  A$, a nonempty  open truncated cone $C_y \subset   \Omega $ with vertex at $y$  such that $\int _{C_y} V(x)\,  \vert  x -y\vert  ^{2-N}\, dx< \infty $.  
\end{theorem}

The next section is mainly devoted to the description of some known facts, important in our approach. In Section 3 some  characterizations of {\sl fine} regularity with respect to a potential $V \in   { \mathcal V}(\Omega)$ are derived. In Sections 4 and 5,  we prove Theorem \ref{th1} while  Section 6 is devoted to the proof of Theorem \ref{th2} and other properties related to almost everywhere regularity. Finally  in Section 7 we observe that  the main results can be   extended to more general pairs of elliptic operators in $\Omega $, ${ \mathcal L}_j={ \mathcal L}- V_j$,  $j=1,\,2$, with ${ \mathcal L}$ in the standard form ${ \mathcal L}(u)= \sum \partial _i(a_{ij}\partial _ju)$ and $0 \leq V_ 1 \leq V_2$. 
\

{\em Acknowledgments.} The author is grateful to Moshe Marcus and Laurent V\'{e}ron for  attracting him  by   motivating questions to the topics considered in this paper.

\section{Boundary Harnack principle and Fatou's  theorem for $\mathbf {L_V}$}

As in the appendix in \cite{VYA} we will use the results of \cite{anc3} in forms which are more or less implicit in \cite{anc3} (see also \cite{ancbook}),  a difference with \cite{VYA} being that we will also draw upon the available Fatou-Doob theorem established in \cite{anc3}. 

For the readers convenience, we would like to state and make clear these (well-known to experts) ancillary results in  forms suitable  for our approach. As a result there will be some overlap  with the exposition in \cite[ Section A.2]{VYA}.

{\bf 2.1.} Let $\;L= \sum_{1 \leq i,j \leq N} \partial _i(a_{ij} \partial _j . ) +  \sum _{j \leq N}b_j\, \partial _j. + \gamma .\;$ be a second order  elliptic operator in divergence form in the open unit ball $B_N$ of ${ \mathbb  R}^N$ and let $ \delta (x)= \delta _{B_N}(x)=1- \vert  x \vert  $ for $x \in   B_N$. We   assume that the coefficients $a_{ij}$, $b_i$ and $ \gamma $ are real measurable in $B_N$ and that for some constant $C>1$ and  all $x \in   B_N$, $ \xi   \in   { \mathbb  R}^N$ : (i)     $   C^{-1} \Vert  \xi   \Vert ^2 \leq \sum _{i,j}   a_{ij}(x)  \xi _i\, \xi  _j \leq  C \Vert  \xi   \Vert ^2$, (ii) $ \delta (x) \sum_j  \vert  b_j(x) \vert  \leq C$ and (iii) $ \delta^2 (x)\,   \vert  \gamma (x) \vert   \leq C$.  It is also assumed that ${ \mathcal L}= \delta ^2\,L$ is weakly coercive, i.e.\,: (iv) there exists $ \varepsilon_0 >0$ such that ${ \mathcal L} + \varepsilon_0 I$ admits a non trivial positive supersolution (ref.  \cite{anc3}). 

Note that (iv) holds in particular if $b\equiv 0$ and $ \gamma  \leq 0$ (assuming (i) and (iii)): thanks to Hardy's inequality (\cite{steinbook}  p. 272 and \cite{kadkuf}) there exists $ \varepsilon = \varepsilon (N,C)>0$ such the form  $ a( \varphi  , \varphi  )=\sum \int_{B_N} a_{ij} \partial _j \varphi  \, \partial _i \varphi  \, dx -\int_{B_N}  (\gamma+ { \frac {   \varepsilon } {  \delta ^2}} )\,  \varphi  ^2 \, dx $, $   \varphi   \in   H_0^1(B_N)$, is coercive in $H_0^1(B_N)$. So if $ \psi  \in   H_0^1(\Omega )$, $ \psi  \geq 0$, Stampacchia's projection theorem \cite{sta0} applied to the convex set $C_ \psi = \{ v \in   H_0^1(\Omega )\,;\, v \geq  \psi \,  \}$ provides an element $u \in   C_ \psi $ such that  $a(u,  \theta  ) \geq 0$ for every $ \theta   \in   [H_0^1(\Omega )]_+$. Thus, $L(u)+ { \frac {   \varepsilon } { \delta ^2 }}u \leq 0$ in $B_N$.

Let $L$ be as above satisfying (i)-(iv) and let $ \Gamma  $ denote the Green's function for $L$ in $\Omega $ -- which exists (see e.g.\ \cite{ancbook}). Thus if $ \Gamma  _y(x):= \Gamma  (x,y)$, we have $L \Gamma  _y=- \delta _y$ in the weak sense \cite{sta}. As noticed in \cite{anc3}, the operator ${ \mathcal L}=(1- \vert  x \vert  ^2)^2 L$ is an ``adapted" elliptic operator in $B_N$ for the hyperbolic metric $g_h(dx)= 4\vert  dx \vert  ^2/(1- \vert  x \vert  ^2)^2$. Since ${ \mathcal L}$ is weakly coercive in $B_N$  one may hence apply the main results in \cite{anc3} (Theorem 3 and Theorem 4) and get the following (notice that ${ \mathcal L}$-Green's function with respect to the hyperbolic volume  is $G_{ \mathcal L} (x,y)= 2^{-N}(1- \vert  y \vert  ^2)^{N-2}\;  \Gamma  (x,y)$): 

(a) for each $P  \in   \partial B_N$, the limit $K_P^L(x)=\lim_{y  \to    P}  \Gamma  _y(x)/ \Gamma  _y(0)$, $x \in   B_N$, exists and defines a positive $L$-harmonic function $K_P^L$  in $B_N$ which depends continuously on $P$ (for uniform convergence on compact subsets of $B_N$), 

(b) Every positive $L$-solution $u$  in $\Omega $ can be written in a unique way as $u(x)=\int _{\partial B_N}\, K_P^L(x)\, d\mu  (P)$ for some positive (finite) measure $\mu  $ in $\partial B_N $,

(c) {\bf Relative Fatou Theorem.} It follows from (b) and the general relative Fatou-Doob-Na\"\i m theorem (see e.g.\ \cite{ancbook} p. 28, Th\'{e}or{\`e}me 1.8) that if $u$, $\mu  $ are  as above and if $s$ is a positive $L$-superharmonic function then $s/u$ admits at $\mu  $-almost every $P \in   \partial B_N$ the  ``fine" limit $f(P)$ (in the sense of the potential theory with respect to $L$, with $P$ seen as a minimal Martin boundary point, see \cite{ancbook}). Here $f={ \frac { d\nu   } {d\mu  }}$ $\mu$-a.e.,  $\nu $ being the measure in $\partial B_N$  associated to the largest $L$-harmonic minorant of $s$. If $s$ is $L$-harmonic, then for every $P \in   \partial B_N$  such that  ${ \frac { s  } {u}} \to  \ell$ finely at  $P \in   \partial B_N$, the ratio  ${ \frac { s(x)  } {u(x)}}$ tends to $  \ell$ as $x \to  P$ nontangentially in $B_N$ \cite{ancbook}.

(d) There is a constant $c=c( \varepsilon_0 ,C,N) \geq 1$ such that for $z=(z',z_N)  \in    B_N $ with $z_N \geq \frac 1 3$ and $y \in   B_N^-:= \{ (x',x_N) \in   B_N\,;\, x_N \leq -{ \frac { 1  } {2}}\,  \}$ one has  \begin{align} \label{phfIni} c^{-1}\,  \Gamma  (z,y) \leq  \Gamma  (z,0)\, \Gamma  (0,y) \leq c\,  \Gamma  (z,y).\end{align}
 This follows from  Theorem 1 in \cite{anc3} and Harnack inequalities, on observing that the distance of $O$ to the hyperbolic geodesic $\buildrel \frown \over  {zy}$ is bounded by a constant.  Letting $z$  tend to a limit position $Q \in   (\partial B_N)^ +:=  \{ x \in   \partial B_N\,;\,, x_N \geq { \frac { 1  } {3}} \}$ it follows that $c^{-1} K_Q(y) \leq   \Gamma  (0,y) \leq c\, K_Q(y)$ for $y \in   B_N^-$. Whence the next proposition.

\begin{proposition} \label{HBhy} Let $u=\int K_Q\, d\mu  (Q)$ where $\mu  $ is a positive measure supported by $(\partial B_N)^+$. Then, for every $y \in   B_N^{-}$ and some constant $c=c( \varepsilon_0 ,C,N) \geq 1$,
\begin{align} c^{-1} u(0)\, \Gamma  (y,0) \leq u(y) \leq c\, u(0)\, \Gamma  (y,0).
\end{align}
In particular, if $u$ and $v$ are positive $L$-harmonic functions in $ B_N$ whose associated measures on $\partial B_N$ are supported by $(\partial B_N)^+$  
\begin{align}\label {eq2.3} \frac {u(x)}{u(0)} \leq \;c^2\; \frac {v(x)}{v(0)}, \;\; x \in   B_N^-.
\end{align}
\end{proposition}

\begin{remark} \label{remark2}{ \rm Suppose $L(\mathbf 1)= \gamma   \leq 0$. Let $u$ be a positive solution of $L(u)=0$ in $B_N$   and let $\mu  $ be the corresponding measure on $\partial B_N$. If $u$ vanishes on $ \{ z \in   \partial B_N\,;\, z_N <{\frac  { 1} {3}}  \}$ then  $\mu  $ is supported by $ \{ z \in   \partial B_N\,;\, z_N \geq {\frac  { 1} {3}}  \}$. This follows  from the Fatou theorem in (c) : ${ \frac { 1  } {u}}$ admits a finite fine  limit at $\mu  $-a.e.\ point $P \in   \partial B_N$.} 
\end{remark}

\vspace{5mm}

{\bf 2.2.} We now apply Proposition \ref{phfIni} to the operators $L_V$ and fix more notations. Given positive reals $r,\,  \rho >0$ such that $0<10\,r< \rho $ and   a ${ \frac {   \rho } {10r }} $--Lipschitz function $f$ in the ball $B_{N-1}(0,r)$ of ${ \mathbb  R}^{N-1}$ such that $f(0)=0$, define :
\begin{equation} U_f(r, \rho ):= \{ (x',x_N) \in   { \mathbb  R}^{N-1}\times { \mathbb  R}\simeq { \mathbb  R}^N\,;\,  \vert  x' \vert  <r,\, f(x')<x_N< \rho \,  \}    \end{equation}
\noindent We will also denote this region $U_f$ (leaving $r$ and  $ \rho $ implicit) or even $U$ when convenient. Set   $\partial _{\#}U:=\partial U \cap  \{ x=(x',x_N) \in   { \mathbb  R}^N\,;\,  \vert  x' \vert  <r,\, x_N=f(x')\,  \}$ and define $T(t):=B_{N-1}(0;tr)\times (-t \rho , +t \rho )$.

The next boundary Harnack principle already mentioned in \cite [Appendix]{VYA} is implicit in \cite{anc3}.

\begin{lemma} \label{PHF}  Let $V \in   { \mathcal V}(U,a)$ and set $L :=  \Delta  -V$.  There is a constant $C$ depending only on $N$, $a$ and ${ \frac {   \rho } {r }}$ such that for  any  two positive $L $-harmonic functions  $u$ and $v$ in $U $ that vanish in $\partial_\# U $, 
 \begin{align} { \hspace {10truemm}  \frac {u(x)} {u(A) }} &\leq  C\; { \frac {  v(x)} {v(A)}} \ \ {f\!or\  all\ \,} x \in   U \cap T({\frac  { 1} {2}}) 
\end{align}
where $A=A_U=(0,\dots, 0,{ \frac {   \rho } {2 }})$.
\end{lemma}

\noindent {   \em Proof.} Let us briefly describe an argument that yields this result starting from  (\ref{phfIni}).
By homogeneity we may assume that $r=1$ and $ \rho $ is fixed.  
Set $A'=(0,...,0,{ \frac { 2  \rho } {3 }})$. 

Using  a  biLipschitz map $F: U \to  B_N(0,1)$ with biLipschitz constants depending only on $ \rho $ and $N$ and mapping $A'$ onto $0$, $U \cap T(1/2)$ onto  $B_N^-$, $\partial U\setminus \partial _\#U$ onto $(\partial B_N)^+$  we are reduced to (\ref{eq2.3}). The function  $u_1:=u \circ F^{-1}$ is ${L}_1-V \circ F^ {-1}$ harmonic for some (symmetric) divergence form elliptic operator $L_1= \sum _{i,j} \partial _i(a_{ij} \partial _j)$ satisfying $ C_1^{-1}\, I_N \leq   \{ a_{ij} \}  \leq C_1\, I_N$, $C_1=C( N,\rho ,a)$. Thus, thanks to remark  \ref{remark2}, Proposition \ref{HBhy} applies to $u_1$, $v_1=v \circ F_1^{-1}$ and $L:=L_1-V_0$, $V_0=V \circ F^{-1}$. The Lemma follows (by Harnack inequalities for  $L$ we may  replace $A$ by $A'$ in Lemma \ref{PHF}). $\square$

\begin{remark}\label {vanish} {\rm Notice that we have also shown that $u(x)/u(A) \leq c \,r^{N-2}\, G_A^V(x)$ in $U \cap T({ \frac { 1  } {3}})$ where $G_A^V$ is Green's function in $U$ with pole at $A$ and with respect to $ \Delta  -V$. Here $c=c(N, a, { \frac { r  } { \rho }})$ is another positive constant.} \end{remark}

\vspace{5mm}

{\bf 2.3.}  By standard arguments (see e.g.\  \cite{anc1}, \cite{anc2}) it follows from Lemma \ref{PHF}    that properties (a) (b) and (c) in {\bf 2.1} above extend to the potential theory in a Lipschitz domain $\Omega $ with respect to $L_V:= \Delta  -V$, $V \in   { \mathcal V}(\Omega ,a)$, $a>0$ -- replacing $B_N$ by $\Omega $. This  is  straight\-forward when $\Omega $ is biLipschitz equivalent to $B_N$ using a biLipschitz  change of variable. It also follows in general by the results in \cite{ancbook} from the Gromov hyperbolicity of $\Omega $ equipped with its pseudo-hyperbolic metric (see \cite[Th.\ 3.6] {BHK} ).

To be more explicit, fix  $x_0 \in   \Omega $ and denote $G^V$ the Green's function for $ \Delta  -V$ in $\Omega $.  Then for $y \in   \partial \Omega $, the limit  $K_y^V(x):=\lim_{z \to  y} G^V(x,z)/G^V(x_0,z)$ exists. Moreover $K_y^V$ depends continuously on $y$, $K_y^V=0$ in $\partial \Omega \setminus  \{ y \}$  and every positive $L_V$-harmonic function $u$ in $\Omega $ can be uniquely written as $u(\cdot)=\int _{\partial \Omega }\, K^V_y(\cdot )\, d\mu  (y)$ where $\mu  $ is a positive finite measure in $\partial \Omega $.

\vspace{5mm}

{\bf 2.4.} The relative Fatou theorem (as stated in {\bf 2.1} (c)) extends as well, but we 
will also need a simple extension  which follows from the same argument: if $u$ and $\mu  $ are as above and if $s:B_N  \to  { \mathbb  R}_+$ satisfies the ``strong Harnack property": 

\begin{center}$(HP)_{\Omega }$ \ \ $\lim_{ \varepsilon \downarrow 0}[\sup  \{ \,s(x)/s(y)\,;\,x,\,y \in   \Omega ,\,  \vert  x-y \vert   \leq  \varepsilon\;  \delta _{\Omega }(x)\, \}]=1$\end{center}

then for $P \in   \partial \Omega $ such that $s/u$ admits a fine limit $\ell$  at $P$, again $s/u$ does admit $\ell$ as  a nontangential limit at $P$.  See e.g.\ the proof of Th\'{e}or{\`e}me 6.5 in \cite{ancbook} p.\ 100. Whence the next proposition.

 \begin{proposition} \label{propRFT}Let  $u$ and $\mu  $ be as above and let $s$ be a nonnegative $L_V$-superharmonic in $\Omega $ satisfying the strong Harnack property $(HP)_{\Omega }$. Then for $\mu  $-almost all $y \in    \partial  \Omega $, the ratio $s/u$ admits a nontangential limit $f(y)$ at $y$, and $f(y)={ \frac {  d\nu } { d\mu  }}(y)$ $\mu  $-a.e., where $\nu $ is the measure on $\partial \Omega $ associated to the largest $L_V$-harmonic minorant of $s$.
\end{proposition}

\vspace{5truemm}
{\bf 2.5.} As a first and direct application, an answer is given to a question  of Moshe Marcus  asking if a positive solution of $ \Delta  u-Vu=0$ that converges nontangentially to zero at every boundary point must be zero.

\begin{proposition} \label{MM} Let $u$ be a nonnegative $L_V$-solution in $\Omega $. Suppose that  for every $P
 \in   \partial \Omega $ and every $ \varepsilon >0$ there is a sequence $ \{ x_n \}$ in $\Omega $ converging nontangentially to $P$ and such that  $\liminf u(x_n) \leq  \varepsilon $. Then $u= 0$ in $\Omega $.
\end{proposition}

{\em Proof.} Since the function $\mathbf 1$ is $L_V$-superharmonic and has the strong  Harnack property in $\Omega $, the ratio ${ \frac { 1  } {u}}$  has a finite nontangential limit at $\mu  $-a.e.\ point $P \in   \partial \Omega $ - if $\mu  $ is the positive measure in $\partial \Omega $ associated to $u$. If $\mu  \ne 0$, this contradicts the assumption on $u$. Thus  $\mu=0  $. $\square$

We may go a little further and show in the same way: if $w$ is a difference of nonnegative $L_V$-harmonic functions in $\Omega $ and if for each $y \in   \partial \Omega $ there is a sequence $ \{ x_n \}$ converging to $y$ nontangentially in $\Omega $  such that $\lim w(x_n)=0$, then $w=0$. To see this write $w=u-v$ with $u=K_\mu ^V $, $v=K_\nu ^V$ and $\mu  \wedge \nu =0$ ($\mu  $ and $\nu $ are positive finite measures in $\partial \Omega $). If $\mu  \ne 0$ then from ${ \frac {  1} {w }}={ \frac {  1} { u}}\, { \frac {  1} {1-{ \frac {  v} {u }} }}$ and Proposition \ref{propRFT}, we get that ${ \frac {  1} {w }}$ has a finite nontangential limit $\mu  $ a.e.\ in $\partial \Omega $ which is absurd. So $\mu  =0$ and similarly $\nu =0$.

\begin{remark} { \rm Let $\omega $ be an open subset of $\partial \Omega $ and let again $u=K_\mu  ^V$ be a positive $L_V$-solution in $\Omega $ such that for every $P \in   \omega $, $u(x)$ admits the inferior limit $0$ as $x \to  P$ nontangentially. Then the argument in Proposition \ref{MM} shows that $\mu  (\omega )=0 $. Moreover $\lim_{x \to  P_0} u(x)=0$ for $P_0 \in   \omega $: assuming as we may  $\Omega =U$, $\omega =\partial U \cap T({ \frac {  1} {4 }})$ and $P_0=0$ with the notations in 2.2, this follows from remark \ref{vanish} since $G_A^V \leq G_A^0$ .}
\end{remark}

\begin{remark}\label{remark8} {\rm In view of Section 7, we note that  the properties  in 2.3 to 2.5  (in particular (a), (b), (c), Prop.\ \ref{HBhy} and \ref{propRFT}) hold  as well for an elliptic operator ${ \mathcal L}_V= \sum _{i,j} \partial _i(a_{ij}\partial _j.) -V$ with $a_{ij}$  measurable in $\Omega $, $a_{ij}=a_{ji}$, $ c^{-1} I \leq \{ a_{ij} \} \leq cI$ in $\Omega $,  $V \in   { \mathcal V}(a,\Omega )$, $c,\, a>0$. The constants in Lemma \ref{PHF} depend then also on $c$.}
\end{remark}

\section{  Boundary behavior of solutions and fine regularity of boundary points.}

 In all this section $\Omega $ is a bounded Lipschitz domain in ${ \mathbb  R}^N$  and $V \in   { \mathcal V}(\Omega ,a)$, $a>0$. 
\vspace{5mm}

\vspace{5mm}

{\bf 3.1.} The next two propositions are  simple consequences of  the available relative Fatou theorem.

\begin{proposition} \label{propA5}If $y_0 \in   \partial \Omega $ is a {\em {finely regular}} boundary point with respect to $L_V$ then  $ \tilde K_y^V(x)/ K_y(x)  \to  1 $ as $x \to  y_0$ non-tangentially. More generally, if $\mu  $ is a  (positive finite) measure  on $\partial \Omega $, then for $\mu  $-almost all $y \in   { \mathcal R}\mbox{\rm \footnotesize{eg}}_{_ V}(\Omega ) \subset   \partial \Omega $ the ratio $ \widetilde K_\mu  ^V(x)/K_\mu (x) $ converges to $1$   as $x$ converges nontangentially  to $y$. 
\end{proposition}

Here, by definition, for any positive finite measure $\mu  $ on $\partial \Omega $,  $ \tilde K_\mu  ^V(x)=\int  \tilde K_y^V(x)\, d\mu  (y)$ and $ K_\mu  ^V(x)=\int   K_y^V(x)\, d\mu  (y)$. Recall that for $y \in   \partial \Omega $, $ \tilde K_y^V$ is the largest $L_V$-subharmonic minorant of $K_y$, that $ \tilde K_y^V=c(y)\,K_y^V$, $0 \leq c(y) \leq 1$, and that $c(y)>0$ if and only if $y \in   { \mathcal R}\mbox{ \footnotesize{eg}}_V(\Omega )$. In particular $ \tilde K_\mu  ^V=K_{c\mu  }^V$.

{\em Proof.} 
The function $p:= K_\mu  -\tilde K_\mu  ^V$ is an $L_V$-potential in $\Omega $ (since $K_y- \tilde K_y^V$ is an $L_V$-potential and potentials are stable under integration) \cite{brelot}. Thus the greatest $L_V$-harmonic minorant of $K_\mu  $ is $ \tilde K_\mu  ^V=K_{c\mu  }^V$. Since $K_\mu  $ satisfies the Harnack property it follows using  Proposition \ref{propRFT} that  $K_\mu (x) / \tilde K_\mu  ^V (x) \to 1$ as $x$ converges nontangentially to $y$, for $\mu  $-almost all $y \in   \partial \Omega $ such that $c(y)>0$. $\square$

\begin{proposition} \label{propA6}  If $\mu  $ is a  (positive finite) measure  on $\partial \Omega $, then for $\mu  $-almost all $y \in   { \mathcal S}\mbox{\rm \footnotesize{ing}}_{_ V}(\Omega ) \subset   \partial \Omega $ the ratio $  K_\mu  ^V(x)/K_\mu (x) $ tends to $+ \infty $ as $x \to  y$ nontangentially in $\Omega$. 
\end{proposition}

{\em Proof.} Set $\mu  '=\mathbf 1_{{ \mathcal S}\mbox{\scriptsize {ing}}_V(\Omega)}\, \mu  $ and assume as we may that $\mu  '\ne 0$. For  $y \in   { \mathcal S}\mbox{\footnotesize{ing}}_V(\Omega )$, $K_y$ is an $L_V$-potential. So (potentials being stable by integration) $K_{\mu '} $ is an $L_V$-potential and by Proposition \ref{propRFT} , ${ \frac {  K_{\mu  '}} { K_\mu  ^V}} \to  0$  nontangentially  at $\mu  $-almost all $y \in   \partial \Omega $.  On the other hand, by Proposition \ref{propRFT} for $V=0$,  ${ \frac {  K_{\mu  }} { K_{\mu '}}}  \to  1$ nontangentially at $\mu  '$-a.e.\ point on $\partial \Omega $. $\square$

Let $\mu  $ be as in Proposition \ref{propA6} and let $u=K_\mu  ^V$. We have essentially shown the following.

 \begin{corollary} \label{coro33} For $\mu  $-almost all $y \in    \partial \Omega $ the ratio $  u(x)/K_\mu (x) $ tends to ${ \frac {  1} {c(y) }} $ as $x \to  y$ nontangentially in $\Omega$. For every positive finite measure $ \theta   $ on $\partial \Omega $, 
the ratio $u(x)/K_ \theta  (x) \to  { \frac {  1} {c(y) }}\,{\frac  { d\mu  } {d \theta   }} (y)$ nontangentially at $y$ for $ \theta   $-almost every $y \in   { \mathcal R}\mbox{\rm \footnotesize{eg}}_V( \Omega) $.
\end{corollary}

{\sl Proof.} (a) The first claim follows from the above and the relation  $u(x)/K_\mu  (x)=$ $ (K_\mu  ^V(x)/ K^V_{c\mu  }(x)) (\tilde K_\mu  ^V(x)/K_\mu  (x))$.

(b) Let $ \theta  $ be a finite positive measure on  $\partial \Omega $. Using the relative Fatou theorem for $ \Delta  -V$ we get that for $ \theta  $-almost all ($V$-finely) regular $y \in   \partial \Omega $:
$$K_{\mu  }^V(x)/K_{ \theta  }(x) =(K_{\mu    }^V(x)/K_{ \theta  }^V(x)) (K_ \theta  ^V(x)/K_ \theta  (x))\to  { \frac {  d\mu  } {d \theta   }}(y)\times {\frac  { 1} {c(y)}} $$ 
as $x \to  y$ nontangentially. $\square$

Corollary \ref{coro33} leads to a natural notion of relative (with respect to $x_0$) {\em fine trace} (compare \cite{dyn1}, \cite{VYA}). Let  $u$ be a positive solution of $ \Delta  u-Vu=0$ in $\Omega $. Let $\mu  $ be the (unique) corresponding positive (finite) measure on $\partial \Omega $ such that $u=K^V_\mu  $. We define the ``fine trace" of $u$ on $\partial \Omega $ as the (nonnecessarily finite) positive  measure $\nu =c^{-1} \mu  $. This measure is $ \sigma  $-finite on the $K_ \sigma  $-set $ \{ c>0 \} $ in $   \partial \Omega $ and  for every nonempty subset of the $G_ \delta $-set $ \{ c=0 \}={ \mathcal S}\mbox{\footnotesize {ing}}_V(\Omega )$ we have (i) $\nu  (A)=+ \infty $ if $\mu  (A)>0$ and $\nu (A)=0$ if $\mu  (A)=0 $. The measure $\nu $ is $ \Sigma $-finite in the sense of \cite[ p.\ 153]{dyn1}. When $\Omega $ is $C^{1, \alpha }$-smooth and the Poisson kernel $P^\Omega $ is well defined, $>0$ and continuous in $\Omega \times \partial \Omega $, it is more natural to define the fine trace of $u$ as $ \tilde \nu = \pi^{-1} \,\nu $, where $  \pi  (y)=P^\Omega (x_0,y)=\partial _{n_y}G_{x_0}$, $n_y$ being the inner normal at $y \in   \partial \Omega $ (for $ \lambda    \in   { \mathcal M}_+(\partial \Omega )$, $P^\Omega _ \lambda   =K( \pi  \lambda   )$). For  results on fine traces related to solutions of semilinear equations, see \cite{dyn1, dyn2, Legall, mselati, MV2, MV3, MV1, MV6, VYA} and references there.

\vspace{5mm}
{\bf 3.2.} Other characterizations of ``fine regularity'' in terms of Green's functions -- already given in the appendix of \cite{VYA} -- are recovered in Proposition \ref{equivdef} below. Before we state a simple variant of the Harnack boundary property Lemma \ref{PHF} and another standard and useful consequence of Lemma \ref{PHF}. Let $U:=U_f(r, \rho )$ be as in {\bf 2.2}. Recall $A=(0,\dots, 0,{ \frac {   \rho } {2 }})$ and $T(t)=B_{N-1}(0,tr)\times (-t \rho ,t \rho )$.

\begin{lemma}\label{eBHI} Let $u$ be  positive harmonic in $U$, let $v$ be positive $ \Delta  -V$-harmonic in $U$, $V \in   { \mathcal V}(a,U)$, and assume  that $u=v=0$ in $\partial_\# U:=\partial U \cap T(1) $. Then 
\begin{align} { \frac {v(x)} {v(A) }} &\leq c\; { \frac {  u(x)} {u(A) }} \hspace {5truemm} \ \ {\rm for\  \ } x \in   U \cap T({\frac  { 1} {2}}) \end{align}
\flushleft and some positive constant $c$ depending only on $ \rho /r$, $N$ and the constant $a$.
\end{lemma}
{\em Proof {\rm (See the Appendix in \cite {VYA})}.}  Assuming as we may $r=1$, we have seen that $v(x) \leq c\,v(A)\, G_{A}^V(x)$ in $U \cap T({ \frac { 1  } {2}})$, $c=c({ \frac {  r} { \rho  }},a,N)$  and we know that $G_A^V \leq G_A^0$ in $U$ (if $G_{A}^V$ is the $( \Delta  -V)$-Green's function in $U$ with pole at $A$). By the maximum principle, the Harnack inequalities and the well-known fact that  $G_A^0(x) \leq c_1:=c_1(r, N)$ in $\partial B(A',{ \frac { r  } {4}})$) we deduce that $u(x) \geq c_1\,v(A)\, G_{A'}^0(x)$ in $U\setminus B(A,{ \frac { r  } {4}})$. So that -- using  Harnack inequalities in $B(A, { \frac { r  } {2}})$ for $u$ and $v$ -- the lemma follows. $\square$

\noindent Denote $g_{x_0}^V $ the Green's function with respect to $ \Delta  -V$ in $\Omega $ with  pole at $x_0$. 

\begin{lemma} \label{GKlemma} Let $y \in   \partial \Omega $ and let $\nu _y \in   { \mathbb  R}^N$ be a pseudo-normal at $y$ for $\Omega $, i.e., $\nu _y\ne 0$ and for some small $  \eta >0$, $C(y,\nu _y,  \eta ):= \{ y+t(\nu _y+v)\, ;\, 0<t \leq   \eta ,\,  \Vert v \Vert  \leq       \eta   \; \} \subset   \Omega\setminus  \{ x_0 \} $. There is a constant $C =C(\Omega ,a, y,\nu _y,   \eta )\geq 1$ such that for $x=y+t\nu _y$, $t \leq    \eta $, 
\begin{align}
C^{-1}\, d(x,y)^{2-N}  \leq K^V_y(x) \;g^V_{x_0}(x) &\leq C\, d(x,y)^{2-N} . \nonumber \end{align}
\end{lemma}

{\em Proof \rm (See \cite{ancbook} p.\ 99, or \cite[Appendix] {VYA}).}  
 We recall the idea of the proof  assuming as we may $  0<\eta  \leq \Vert \nu _y \Vert=1$ and $ \vert  x_0-y \vert   \geq 2   \eta   $. Let  $u=K^V_y $ and $v(x):=g^V_{y+t\nu_y }(x)$ for a given $t \in   (0,  \eta ] $. Using  the estimate $v\sim t^{2-N}$ in $\partial B(y+t\nu_y ,  {\frac  {  t   \eta } {2}}  )$,     one  sees that $u(x) \sim u(y+t\nu_y )\, t^{N-2}\, g_{y+t\nu_y }^V(x)$ for $x \in   \Omega  \cap  \partial B(y, t ( \Vert\nu _y \Vert + \eta /2   ))$ (here $\sim$ means ``is in between two constant times" with constants depending only on $y$, $\Omega $, $\nu_y $, $  \eta $  and $a$). By the boundary Harnack principle (lemma \ref{PHF}) it follows that $$u(x) \sim u(y+t\nu_y )\, t^{N-2}\, g_{y+t\nu_y }^V(x)$$ for $x \in   \Omega  \setminus  B(y, t ( \Vert\nu _y \Vert + \eta /2   ))$. Taking  $x=x_0$, it follows that $K^V_y(y+t\nu_y )\sim  1  /(t^{N-2}g^V(y+t\nu_y ;x_0)) $. $\square$ 

\begin{corollary} \label{KVKcoro}Under the  assumptions of Lemma \ref{GKlemma} there is a constant $C'>0$ such that for all small $t>0$,
\begin{align} K_y(y+t\nu _y) \leq C'\, K_y^V(y+t\nu _y).\end{align} 
\end{corollary}

\begin{proposition} \label{equivdef} {\rm ({\rm \cite[Appendix]{VYA}})} Given $y \in   \partial \Omega$ and a pseudo-normal $\nu _y$ at $y$ for $\Omega $, the following are equivalent:
\begin{enumerate}
\item  \  $\tilde K^V_y=0$ (i.e.\, $y \in   { \mathcal S}\mbox{\rm \footnotesize{ing}}_V(\Omega )$),

\item \ $\limsup_{t\downarrow 0} K_y^V(y+t\nu _y)/K_y(y+t\nu _y)=+ \infty $,

\item $\lim_{t\downarrow 0} K_y^V(y+t\nu _y)/K_y(y+t\nu _y)=+ \infty $,

\item $\lim_{\Omega \ni x \to  y} g_{x_0}^V(x)/g^0_{x_0}(x)=0$.
\end{enumerate}
\end{proposition}

\vspace{3truemm}

{\sl Proof of Proposition \ref{equivdef}.} (a) The equivalence of (i), (ii) and (iii) is already established (using Propositions \ref{propA5} and \ref{propA6}). 
 Using Lemma \ref{GKlemma} we see -- as in the appendix of \cite{VYA} -- that (ii) is equivalent to (iv)$'$: $\liminf_{t\downarrow 0} g_{x_0}^V(y+t\nu _y)/g_{x_0}^0(y+t\nu _y)=0 $. 

(b) To see that (iv)$'$ $\Rightarrow$ (iv) we may assume that $y=0$ and that (with the notations above in {\bf 2.2 }) $T(1) \cap \Omega =U$, $U=U_f(r, \rho )$ with $x_0 \in   \Omega \setminus \overline  U$. Applying Lemma \ref{eBHI} to $U$,  $u=g^V_{x_0}$,$v=g_{x_0}$, and $U_t=U_{t_j}$ for a sequence $t_j$ such that $t_j\downarrow 0$ and $u(A_{t_j})=o(v(A_{t_j}))$, $A_{t_j}=(0,\dots, 0,t_j)$, the desired implication is obtained.  Thus, using (a) again, (ii), (iii) and (iv) are equivalent.
And Proposition \ref{equivdef} is proved.   $\square$

\section{Characterization of regularity}

 Let again $\Omega $ be a bounded Lipschitz domain in ${ \mathbb  R}^N$ and let $V \in   { \mathcal V}_a(\Omega )$, $a>0$. We will prove several characterizations of fine regularity of a boundary point $y$ with respect to $ \Delta  -V$.

Let $G$ denote the (classical) Green's function in $\Omega $. So by definition $G(VK_y^V)(x_0)=\int_\Omega  G(x_0,x)\, V(x)\, K_y^V(x)\, dx$, and similarly  $G(VK_y)(x_0)=\int_\Omega  G(x_0,x)\, V(x)\, K_y(x)\, dx$. 

\begin{remark} {\rm {\footnotesize    The finiteness of $G(VK_y^V)(x_0)=\int_\Omega  G(x_0,x)\, V(x)\, K_y^V(x)\, dx$ (or of $G(VK_y)(x_0)$) depends   only on the behavior of $V(x)$ as $x \to  y$. To see this it suffices to show  that given  $x_1 \in   \Omega $, $x_1\ne x_0$, and  a subdomain $\omega  \subset   \Omega $, such that $y, \,x_0,\, x_1 \notin \overline  \omega $, we have $G({\mathbf 1}_\omega VK^V_y)(x_0)< \infty $ and $G({\mathbf 1}_\omega VK_y)(x_0)< \infty $. By the Harnack boundary principle Lemma \ref{PHF}, $K_y^V(x) \leq C\, G^V(x_1,x)$ in $\omega $ (where $C=C_{\Omega ,\omega , V,x_1}>0$) and 
$$\int_\omega   G(x_0,x)\, V(x)\, K_y^V(x)\, dx  \leq  C\, \int_\Omega  G(x_0,x) V(x) G^V(x,x_1)\, dx  \leq C\,G(x_0,x_1)< \infty $$
thanks to the resolvent equation $G=G^V+GM_VG^V$ where $M_V$ is the ``multiplication by $V$" kernel (over $\Omega $) and where $GM_VG^V$ is the composed kernel $G \circ M_V \circ G^V$.

To prove that $G({\mathbf 1}_\omega VK_y)(x_0)< \infty $, it suffices as before  to show  that $G(VG_{x_1})(x_0)< \infty $ (which is stronger than $G(VG^V_{x_1})(x_0)< \infty $). By a   standard estimate: $G(x,x_1) \leq c \, \vert  \delta _\Omega (x) \vert  ^ \alpha $, for $0< \delta _\Omega (x) \leq  {\frac  { 1} {2}}  \delta _\Omega (x_1)$ and constants $c=c(\Omega ,x_1)>0$, $ \alpha  \in   (0,1)$. Thus $G_{x_1} (x)\, V(x) \leq C\,  \delta (x)^ {\alpha -2}$ in $\Omega $ (away from $x_1$). However, by \cite[p.\ \nolinebreak 286-287]{anc3'},  there is a   smooth positive superharmonic function $s$ in $\Omega $ such that $ \Delta  s(x) \leq -{\frac  {  1} { \delta _\Omega ^{2- \alpha }(x)}}  $ in $\Omega $. This implies that  $G({\frac  { 1} { \,\,\delta _\Omega ^{2- \alpha } }})  \leq s< \infty $  in $\Omega $.}} \end{remark}

\begin{theorem}\label{charac} Let $y \in   \partial \Omega $. The following properties are equivalent:
\begin{enumerate}
\item $y$ is finely regular for $L_V= \Delta  -V$ in $\Omega $,
\item  $G(VK_y^V)(x_0)< \infty $,
\item $K_y^V$ admits a $ \Delta  $-superharmonic majorant in $\Omega $,
\item $G(VK_y)(x_0)< \infty $.
\end{enumerate}
\end{theorem}

\begin{remark} \label{remarkApp}{\rm  We recover the necessary condition for regularity of \cite{VYA} (Appendix)  alluded to in the introduction: if there  is an open truncated cone $C$ with vertex at $y$ such that  $C_  \eta :=  \bigcup_{z \in   C}    B(z;   \eta  \vert  z-y \vert  ) \subset   \Omega $ for some $  \eta  >0$ and $\int _C\, { \frac {  V(x)} { \vert  x-y \vert  ^{N-2} }} dx =+ \infty $, then $y \in  { \mathcal S}\mbox{\footnotesize {ing}}_V(\Omega )$. 
For, on $ C$ and away from $x_0$, we have  
$$K_y^V(z)\,G(x_0,z) \geq { \frac {  C} { \vert  x_0-z \vert  ^{N-2} }}{ \frac {  G(x_0,z)} {G^V(x_0,z) }} \geq  { \frac {  C} { \vert  x_0-z \vert  ^{N-2} }}$$  by Lemma \ref{GKlemma}.  Thus $G(VK_y^V)(x_0) \geq \int _{ C } dx/ \vert  x-y \vert  ^{N-2} $.} \end{remark}
For $\Omega $  sufficiently smooth $(i)\Leftrightarrow(iv)$ yields the following rather explicit criteria.

\begin{corollary} If $\Omega $ is $C^{1, \alpha }$, $0< \alpha  \leq 1$, in a neighborhood of $y \in   \partial \Omega $ then $y \in  { \mathcal S}\mbox{\rm\footnotesize {ing}}_V(\Omega )$ if and only if 
$$ \int _\Omega  { \frac {   \vert     \delta _\Omega (x) \vert  ^2} { \vert  x-y \vert  ^N }}\, V(x)\, dx < \infty .$$
\end{corollary}

This is because   $G(.,x_0)\sim  \delta _\Omega (.)$ and $K_y(.)\sim   \delta_\Omega  (.) \vert  .-y \vert  ^{-N}$ near $y$ in $\Omega $ (where $A\sim B$ means ``$A$ is in between two positive constants times $B$'').

{\bf Proof of Theorem \ref{charac}.}  We first show that $(i)\Rightarrow(ii)\Rightarrow(iii)\Rightarrow(i)$ and $(iv)\Rightarrow (i)$. 

(i) $\Rightarrow$ (ii) If $y$ is finely regular then  for some small $c>0$, $w:=K_y-cK_y^V$ is positive $ \Delta  $-superharmonic and $ \Delta  w=-cVK_y^V$. Whence $G(VK_y^V) \leq w$ in $\Omega $ (by the Riesz decomposition) and $G(VK_y^V)(x_0)< \infty $. 

(ii) $\Rightarrow$ (iii): By assumption,  $p:=G(VK_y^V)$ is a $ \Delta  $-potential (i.e., $p\not \equiv + \infty $ in $\Omega $). Thus, the function $w:=K_y^V+p=K_y^V+G(VK_y^V )$ is $ \Delta  $-harmonic, since $ \Delta  w=VK_y^V-VK_y^V=0$.  Since $w \geq K_y^V $ in $\Omega $, this proves  (iii).

(iii) $\Rightarrow$ (i) The function $K_y^V$ is subharmonic in $\Omega $ and by the assumption there is a   smallest $ \Delta  $-superharmonic majorant of $K_y^V$ in $\Omega $ which we denote $u$. The function  $u$ is $ \Delta  $-harmonic and $u=K_y^V+p$, where $p$ is a $ \Delta  $-potential on $\Omega $. Since $u$ is positive harmonic, $u(\cdot)=\int  K_z(\cdot)\, d\mu  (z)$  for some positive finite measure $\mu  $ on $\partial \Omega $ and 
\begin{equation}\label{eq}\int  K_z\, d\mu  (z)   \geq  K_y^V 
\end{equation}
in $\Omega $. Viewing  each member of (\ref{eq}) as an integral of nonnegative $ (\Delta  -V)$-super-harmonic functions  we get that $ \int \nu _z\, d\mu  (z)  \geq  \delta _y$ where $\nu _z$ is the restriction to $\partial \Omega $ of the measure  in $\overline \Omega $ associated to the Green-Martin decomposition of $K_z$ (with respect to $ \Delta  -V$). However, $\nu _z(\partial \Omega \setminus  \{ z \})=0$ ($K_z=c_zK_z^V+(K_z-c_zK_z^V)$ is the Riesz decomposition of $K_z$ with respect to $ \Delta  -V
$ in $\Omega $). Thus $\mu  $ and $\nu _y$ must charge $y$. In particular $\nu _y=c\,  \delta _y$ for some real $c>0$ and  $ K_y  \geq c\,K_y^V$ as desired.

Thus $(i)\Leftrightarrow (ii) \Leftrightarrow (iii) $.

Next we show that $(iv)\Rightarrow (i)$. Let $p=G(VK_y)$. The function $p$ is by  assumption $(iv)$ a $ \Delta  $-potential in $\Omega $. Moreover 
\begin{equation} \Delta  (K_y-p)-V(K_y-p)=VK_y-V(K_y-p)=Vp \geq 0.
\end{equation} Thus $K_y-p$ is $( \Delta  -V)$-subharmonic. Since  $K_y$ is a $ (\Delta  -V)$-superharmonic majorant of this function, it follows that there is a function $u$ which is $ (\Delta  -V)$-harmonic in $\Omega $ and such that 
 $K_y-p \leq u \leq K_y$. Thus the greatest harmonic minorant $ \tilde K_y^V$ of $K_y$ with respect to $( \Delta  -V)$  is larger than $ K_y-p$ in $\Omega $ and hence $ \tilde K_y^V\not\equiv 0$ for otherwise $K_y \leq p$ which is absurd since $p$ is a $ \Delta  $-potential.  Thus $(i)$ holds.

To finish the proof we must show that $(i)$ (or $(ii)$ or $(iii)$) implies $(iv)$. This is the aim of the next section.

\section{Proof of the ``$\mathbf{ G(VK_y)< +\infty }$" criterion}

The proof that $(i)$  implies $(iv)$ (as stated in Theorem \ref{charac})  relies on an  estimate given in Proposition  \ref{propmain} below. A main tool is the following Lemma \ref{mainlemma} (and its corollary) which is interesting on its own right.

We resume to  the situation considered  in Section {2.2} :   $r,\,  \rho >0$ are  such that  $0<10\,r< \rho $,  $f$ is  ${ \frac {   \rho } {10r }} $-Lipschitz in $B_{N-1}(0,r)$  such that $f(0)=0$  and
\begin{equation} U=U_f(r, \rho ):= \{ (x',x_N) \in   { \mathbb  R}^{N-1}\times { \mathbb  R}\simeq { \mathbb  R}^N\,;\,  \vert  x' \vert  <r,\, f(x')<x_N< \rho \,  \}. \nonumber    \end{equation}
\noindent Recall that  for $0<t \leq 1$, $T(t):=B_{N-1}(0;tr)\times (-t \rho , +t \rho )$, $U_t=U \cap T(t)$, and that  $\partial _{\#}U:=\partial U \cap T(1)$.

\begin{lemma} \label{mainlemma}Let $0<t_1 \leq  1$, let $ \zeta_0  \in   W:=U _{t_1}$ and let 
$w(x)=G(x, \zeta_0  )$ where $G$ is the standard Green's function in $W:=U_{t_1}$. Then for every   measurable subset $A$ of $W$,
 \begin{align} \label{Eq5.1} \int _A ({\frac  { w(x)} { \delta _{W} (x)}})^2 \, dx  \leq C\,\; ^{W}\!  \widehat R_w^A( \zeta_0  ) 
\end{align}  where $C$ is a positive constant depending only on  ${ \frac {   \rho } {r }} $ and $N  $ (and not on the set $A$).
\end{lemma}

The r\'{e}duite $^{W}\!R_w^A $  is by definition the lower envelope of all positive superharmonic functions $s$ in $W$ such that $s \geq w$ in $A$. It is harmonic in $W\setminus \overline A$  and its l.s.c.\  regularization $^{W}\! \widehat R_w^A $  is superharmonic in $W$, in fact the smallest nonnegative superharmonic  function in $W$ which is larger than $ w$ in $A$ except at most on a polar subset. Ref.\ \cite{brecl}, \cite{brelot}, \cite{doob}. 

Proposition \ref{mainproposition2} and  inequality (\ref{eq71}) thereafter provide two extensions of Lemma \ref{mainlemma}.

{\bf Proof.} We may  assume $t_1=1$ and $W=U$. Standard  approximation arguments show that we may also assume that $A$ is open and  relatively compact in $W\setminus  \{  \zeta  _0 \}$. Then $ {^W}\! \widehat R_w^A=  {^W}\!R_w^A$ and ${^W}\!R_w^A$ is a Green potential in $W$.

We use an idea from \cite[section VII]{ancbook0} which exploits an identity in the Dirichlet space $H_0^1(U )$ and a familiar Hardy's inequality, $H_0^1(U )$ being equipped with the Dirichlet scalar product $ \langle  \varphi  , \psi  \rangle  _H=\int_U  \, \nabla  \varphi  .\nabla  \psi \, dx$.  Let $ \varphi_0  $ denote a   nonnegative   function in $C_0^ \infty ({ \mathbb  R}^ N)$ with support in $B( \zeta  _0;  \varepsilon )$, $0< 4\varepsilon <d( \zeta  _0,A \cup \partial W)$  and such that $\int  \varphi_0  \, dx =1$.  Let $p$ denotes  its Green's potential in $W $, i.e.\, $p(x)=G( \varphi_0  )(x)=\int G(x,y)\,  \varphi_0  (y)\, dy$. Note that, by Harnack inequalities and the maximum principle, $p(x) \geq { \frac {  1} {c }}G(x,  \zeta  _0)$ for $ \vert  x- \zeta  _0 \vert   \geq 2 \varepsilon $ and some $c=c(N) \geq 1$.

As well-known from Dirichlet and Sobolev spaces  theory, $p$ is in $H_0^1(W ) =   H_0^1(U)$ and $s:=  R_{p}^A$ can be seen as the element in  $H:=H_0^1(W )$ such that $s \geq p$  a.e.\ in $A$  with the smallest Dirichlet norm. Moreover,  using Harnack inequalities, we have
\begin{align} \label{eqprel1}  \langle s, G( \varphi  ) \rangle  _H=\int_W s\, \varphi\, dx  \leq c\, s( \zeta_0 ) 
   \end{align}
for a constant $c=c'(N)$. On the other hand
\begin{align}  \langle s, G( \varphi  ) \rangle  _H=   \int_W G( \varphi  )\, d(- \Delta  s) = \int_W s\, d(- \Delta  s)= \langle s,s \rangle  _H
  & \nonumber \end{align}
because $G( \varphi  )=s$ at every point  of $\overline  A \cap U  $ where $A$ is unthin and the set  of these points  supports the measure $- \Delta  s$ which is the swept-out measure of $ \varphi(x)\,  dx$ on $A$ in $U $.  A simpler argument, not using the fine potential theory,  is obtained on approximating $s $ by $s_n=R_s^{A_n}$ where $ \{ A_n \}_{n \geq 1}$ is an increasing sequence of relatively compact open  subsets of $A$ such that $ \bigcup    A_n=A$. For then $ \langle s_n, G( \varphi  ) \rangle  _H=   \int G( \varphi  )\, d(- \Delta  s_n) = \int s\, d(- \Delta  s_n)$ because $s= G( \varphi  )$ in $A$. Letting $n \to   \infty $ in the equality $ \langle s_n, G( \varphi  ) \rangle  _H=\int s_n\, d(- \Delta  s)$ the relation $ \langle s, G( \varphi  ) \rangle  _H=\int s\, d(- \Delta  s)$ follows since $s=\sup _n s_n$.

Using now a well-known generalization to Lipschitz domains (in our case $U $) of an inequality due to Hardy for ${ \mathbb  R}_+$ (see e.g.\ \cite{kadkuf})   we obtain that, for a constant $C=C({ \frac {  r} { \rho  }})>0$,
\begin{align} \label{eqprel2} \langle s,G( \varphi  ) \rangle  _H&= \int_U \,  \vert  \nabla s(x)\vert  ^2 \, dx \geq C\, \int _{U } \, { \frac {   \vert  s(x) \vert  ^2} {  \delta_U  (x)^2}} \, dx  \geq C \, \int _{A } \, { \frac {   \vert  s(x) \vert  ^2} {  \delta_U(x)^2}} \, dx  \nonumber \\ & \geq C\, c^{-1} \, \int _{A } \, { \frac {   \vert  w(x) \vert  ^2} {    \delta_U  (x)^2}} \, dx     \end{align}

since $s(x)   = p(x) \geq { \frac {  1} {c }} w(x)$ in $A$. Putting together (\ref{eqprel1}) and (\ref{eqprel2}) the result follows. $\square$

\begin{corollary} \label{maincoro}Let $0<t<t' \leq 1$ and $ \zeta  _1=(0,\dots, 0,{ \frac {  t+t'} {2 }} \rho )$. If $u$ is an arbitrary positive harmonic function in $W=U_{t'}$ that vanishes on $\partial _\#U \cap \overline W$,  then for every  Borel subset $A$ of $U_t$,
\begin{align}  \int _A ({\frac  { u(x)} { \delta _{W} (x)}})^2 \, dx  \leq C\, r^{N-2}\,u( \zeta  _1)\; ^{W}\! \widehat R_u^A( \zeta_1  ) 
\end{align}  where $C$ is a positive constant depending only on $t'/t$, ${ \frac {   \rho } {r }} $ and $N  $, but not on the set $A$.
\end{corollary}

This follows from Lemma \ref{mainlemma} since, assuming as we may that $t>{ \frac {  t'} {2 }}={ \frac {  1} {2 }}$, the boundary Harnack principle and the known behavior near $ \zeta  _1$ of $w:=g_{ \zeta  _1}$ -- the Green's function with pole at $ \zeta  _1$  in $W$ --  imply that $u(x)$ is equivalent in size to $u( \zeta  _1)\,((t'-t)r)^{N-2}\, w(x)$ in $U_t$,  with equivalence constants depending only on $N$, $t'$  and $ \rho /r$. 

Corollary \ref{maincoro} (with its proof) may be easily extended as follows.

\begin{corollary}\label{maincorogene} {\rm Let $\omega $, $\omega '$ be open subsets of ${ \mathbb  R}^d$ such that $\overline \omega ' \subset   \omega $ and let $x_0 \in   \Omega \setminus \overline \omega $. There is a constant $C=C(\Omega,\omega ,\omega ', x_0 )$ such that for every positive harmonic function $h$ in $\Omega $ vanishing on $\omega  \cap \partial \Omega $ and every measurable  $A \subset   \omega' $
\begin{align}  
\int _A\,({ \frac {  h(x)} { \delta _\Omega (x) }})^2\, dx  \leq C\,h(x_0)\,  ^{\Omega }\!\widehat R_h^A(x_0).
\end{align}}
\end{corollary}

\vspace{5truemm}
 
Using Corollary \ref{maincoro} we  will prove  the next Proposition  \ref{propmain} which is the key for the proof that $G(VK_y)< \infty$ implies  $y \in    { \mathcal R}\mbox{\rm \footnotesize{eg}}_V(\Omega )$. Notations and assumptions are as above. In particular $0<t<t' \leq 1$ and $ \zeta  _1=(0,\dots, 0,{ \frac {  t+t'} {2 }} \rho )$. 
We are given moreover an element  $V $ of $ { \mathcal V}_a(U )$ (for some $a>0$), a positive solution $u$ of $ \Delta  u-Vu=0$ in $W:=U_{t'}$ that vanishes on $ \partial _\#U \cap \overline W$ and a positive solution $v$ of $ \Delta  v=0$ in $ W$ that also vanishes on $ \partial _\#U \cap \overline W$.

\begin{proposition} \label{propmain}    There is a   positive constant $C=C(N,{ \frac {  t'} {t }},{ \frac {   \rho } {r }},a)$ $>0$ such that  
\begin{align}\label{keyeq}  u(  \zeta  _1)\int_{U_t} \,  \vert  v(x) \vert  ^2\, V(x)\,  dx  \leq Cv( \zeta  _1)\, \int _{ W} &v(x)\,V(x)\, u(x)\, dx\nonumber
\\ &= Cv( \zeta  _1) \int _{ W} v(x)  \,\Delta  u(x)\,\, dx  \end{align}
\end{proposition}

More precisely, given $N$, $a$ and $ \rho/r$ the constant $C$ stays bounded for $t \leq  \beta t'$ and $ \beta <1$ fixed.

{\bf Proof.} a) We first make some simple reductions. 

Observe that by enlarging $t$ we may assume that $t>t'/2$.

Set $  W_1=U_{(t+2t')/3}$, $B_1:=B( \zeta  _1,{\frac  { (t'-t)r} {20}} )$, $   \tilde W_1=U_{(t+2t')/3}\setminus B_1$.  It suffices of course to prove  (\ref{keyeq}) with $W$  in the right-hand side (as the integration domain) replaced by $   \tilde W_1=U_{(t+2t')/3}\setminus B_1$,  $v$ and $u$  being kept as before.

Now if $ \tilde V=\mathbf 1_{W\setminus B_1} V$ and if $ \tilde u$ is  a positive $ \Delta  - \tilde V$-harmonic function in $W=U_{t'}$ vanishing on  $ \partial _\#U \cap \overline W$, we have  by the boundary Harnack principle Lemma \ref{PHF} (and Harnack inequalities), an equivalence in size $ c^{-1} {\frac  {  \tilde u( \zeta  _1)} {u( \zeta  _1)}} u(x)  \leq \tilde u(x) \leq c{\frac  {  \tilde u( \zeta  _1)} {u( \zeta  _1)}} u(x)$ in $W_1 \cap     \{ x=(x',x_N)\,;\, x_N \leq t \rho   \}$ with a constant $c=c(N,a,t'/t,r/ \rho ) \geq 1$. By Harnack inequalities for the operators $ \Delta  -V$ and $ \Delta  - \tilde V$, this equivalence extends to $W_1$.

Thus to prove Proposition \ref{propmain} it suffices to deal with the case where  $V$ is supported by $W\setminus B_1$ -- which we assume from now on -- and prove 
\begin{align} \label{mainest} u(  \zeta  _1)\int_{U_t} \,  \vert  v(x) \vert  ^2\, V(x)\,  dx  \leq Cv( \zeta  _1)\, \int _{ W_1} &v(x)\,V(x)\, u(x)\, dx.
 \end{align}

We also assume as we may  that $v( \zeta   _1)=u( \zeta  _1)=1$.

\vspace{4truemm}

b) Let $w$ denote the (positive) harmonic function in $ W_1=U_{ \frac {  t+2t'} {3 }}$ such that $w=u$ in $\partial   W_1$. By the (extended) boundary Harnack principle Lemma \ref{eBHI} we have $u \leq cv$ in 
$\overline  {W_1}$ for some constant  $c=c(N,a,r/ \rho, t'/t )$,  and hence, by the maximum principle, $w \leq  cv$ in $   W_1$. By Harnack inequalities $u$ is of the order of a constant in $\partial   W_1 \cap     \{ x_N= \rho (t+2t')/3 \}$ and clearly $c^{-1} \leq w( \zeta  _1) \leq c$ with another constant $c$.
Notice finally that by the boundary Harnack property (for $V=0$)  we have an equivalence $c^{-1} w \leq v \leq cw$ in $U_t$, with a constant of the same type as before.

c) Denote  $s$ the difference function,   $s(x)=w(x)-u(x)$, $x \in   \overline  {   W_1}$. Clearly $s=  G(Vu)= G^V(Vw)$ where $  G$ (resp.\ $  G^V$) is Green's function for $    W_1$ with respect to $ \Delta  $ (resp. with respect to $ \Delta  -V.$).

Now $s= G(Vu)$ means in particular  that  
\begin{align}
 s( \zeta  _1)  =  \int _{   W_1 } \,   G( \zeta  _1,x)\,\, V(x)\, u(x)\, dx \nonumber \end{align}
and 
$$ s( \zeta  _1)\leq C\, r^{2-N}\,  \int_{ \tilde W_1} v(x)\, V(x) \,u(x)\, dx$$

using the inequality $G( \zeta  _1,x) \leq C\, r^{2-N}\, {\frac  { v(x)} {v( \zeta  _1)}}$ in $ \tilde W_1=W_1\setminus B_1$ which follows from the maximum principle, and the fact that  $V$ is  supported  by $W\setminus B_1$. Here $C=C(N,{ \frac {  t'} {t }},{ \frac {   \rho } {r }})$.

Next consider the set $A= \{ x \in   U_t\,;\, u(x) \leq {\frac  { 1} {2}} w(x)\,  \}$. 

We have $s \geq {\frac  { 1} {2}} w$ in $A$. So $\, ^{W_1} \! R_w^A \leq 2s$ in $W_1$. By the above and  Corollary \ref{maincoro}
\begin{equation} \int _ A ({\frac  { w(x)} { \delta_W  (x)}})^2 \, dx  \leq 2 \,c\,r^{N-2}\,s( \zeta  _1)  \leq 2\, c'\,\int _{W_1} v(x)\, V(x)\,u(x)\, dx. \nonumber
\end{equation}
with $c'=c(N,{ \frac {  r} { \rho  }},{ \frac {  t'} {t }})$. Thus, since $V \in   { \mathcal V}(a, U)$,  
\begin{align}\int_A w(x)^2\, V(x)\, dx \leq a \int _ A ({\frac  { w(x)} { \delta_W (x)}})^2 \, dx\leq 2\,a\, c'\,\int _{ W_1}\, v(x)\, V(x)\, u(x)\, dx.\end{align} 

On the other hand, since on $U_t\setminus A$ we have by  definition $w \leq 2u$,
 \begin{align} \int _{U_t\setminus A} (w(x))^2 \, V(x)\, dx&= \int  _{U_t\setminus A} w(x)\, ({\frac  { w(x)} {u(x)}} )\;  \Delta  u(x)\, dx \nonumber
\\ & \leq 2\, \int _{U_t\setminus A} \, w(x)\,  \Delta  u(x)\, dx \nonumber
\\ & \leq  2\int _{W_1} w(x)\,  \Delta  u(x)\, dx=2\int _{W_1} w(x)\,V(x)\,u(x)\,dx
\end{align}  

and (\ref{mainest}) follows -- recall that $c^{-1} w \leq v \leq cw$ in $U_t$ and $w \leq cv$ in $W_1$. $\square$

{\em End of the proof of Theorem \ref{charac}.} We now show that if $y_0 \in   \partial \Omega $  is a $V$-finely  regular boundary point for $\Omega $ then $G(VK_{y_0})(x_0)< \infty $. Here $G$ is the classical Green's function in $\Omega $ and as before $x_0$ is a fixed point in $\Omega $, $K_{y_0}$ is the Martin kernel in $\Omega $  with pole at $y_0$. We will denote $G^V$ the Green's function for $ \Delta  -V$ in $\Omega $. 

Again we assume without loss of generality  that $y_0=0$, and that for some $ \rho ,\, r>0$, $10\, r< \rho $,  $\Omega  \cap  \{ (x',x_N) \in   { \mathbb  R}^N\,;\,  \vert  x' \vert  <r,\,  \vert  x_N< \rho \,  \}= U_{f}(r, \rho )$ for a ${ \frac {   \rho } {10 r }}$-Lipschitz function $f:{ \mathbb  R}^{N-1} \to  { \mathbb  R}$ such that $f(0)=0$ (see Section {2.2}). The regions $T(t)$ and $U_t$ are defined as before. We note $\nu =(0,\dots, 0,1)$ the natural pseudo-normal at $y_0$.

We  know that $G(VK_0^V)(x_0)< \infty $ and  $K_{y_0}^V \leq  c_1\, K_{y_0}$. Moreover the functions $K_{y_0}$, $K_{y_0}^V$ satisfy 
\begin{equation} \label{K0KV} c^{-1} K_{y_0} \leq K_{y_0}^V \leq c_1\, K_{y_0}
\end{equation}
on the normal ray $ \{ y_0+t \rho \nu \,;\, 0<t \leq { \frac {  1} {2 }} \}$ for some constant $c=c(N,r/ \rho ,a) \geq 1$  -- see Proposition \ref{equivdef}. Note that the first estimate in (\ref{K0KV}) always holds, independently of the fine regularity assumption for $y_0$, as follows from Lemma \ref{GKlemma}.

 Using a Whitney decomposition of the region  $ \{ x' \in   { \mathbb  R}^{N-1}\,;\, 0< \vert  x' \vert  <r \}$ in ${ \mathbb  R}^{N-1}$, we can cover the punctured ball $B'= \{ x' \in   { \mathbb  R}^{N-1}\,;\, 0< \vert  x' \vert   \leq { \frac {  r} {10 }}\, \}$ in ${ \mathbb  R}^{N-1}$ by a family of balls $B_ \alpha :=B_ {N-1}(x'_ \alpha , r_ \alpha ) $, $ \alpha  \in   J$, of ${ \mathbb  R}^{N-1}$, such that : (i)  $0<3 r_ \alpha  \leq \vert  x'_ \alpha  \vert   \leq C\, r_ \alpha $ for each $ \alpha \in   J $, (ii)   $\# \{  \alpha  \in   J \,;\, B_ \alpha \ni x' \,  \} \leq C$ where $C$ is a constant $ \geq 3$ that depends only on $N$ when $x'$ varies in ${ \mathbb R} ^{N-1}$ (or $ \alpha  \in   J$), and (iii) $B_ \alpha  \cap    B'\ne  \emptyset  $, $ \forall    \alpha  \in   J$.

Fix $ \alpha  \in   J$, denote $B'_ \alpha =B(x'_ \alpha , 2r_ \alpha )$ and consider the cylinders 
\begin{align} T'_ \alpha :=B'_\alpha \times [f(x'_ \alpha )-{ \frac {  3 \rho } {2r }} r_ \alpha ,f(x'_ \alpha )+\,{ \frac {   3\rho } {2r }} r_ \alpha ], \;\; T _ \alpha :=B_\alpha \times [f(x'_ \alpha )-{ \frac {   \rho } {r }} r_ \alpha ,f(x'_ \alpha )+\,{ \frac {   \rho } {r }} r_ \alpha ]\nonumber \end{align} which are contained in $T(1)$. Set $U_ \alpha =T'_ \alpha  \cap \Omega =T'_ \alpha  \cap U$, $V_ \alpha =T_ \alpha  \cap U$ and $A_ \alpha =(x'_ \alpha ,f(x'_ \alpha )+{ \frac {  2 \rho } {r }}r_ \alpha )$.

By the boundary Harnack principle  Lemma \ref{PHF}, (\ref{K0KV}) and Harnack inequalities we have for $x \in   U_ \alpha $, 
$$K_ 0^V (x) \geq { \frac {  1} {c }}\, r_ \alpha ^{N-2} K_0^V(A_ \alpha )\,G^V(A_ \alpha , x) \geq { \frac {  1} {c' }}\, r_ \alpha ^{N-2}\,K_0(A_ \alpha )\, G^V(A_ \alpha ,x)$$
in $U_ \alpha $  (an argument involving the maximum principle suffices for the first inequality). Here and below $c$, $c'$ denote positive constants that do not depend on $ \alpha $ and whose value may vary from line to line.

Similarly $G(x_0,x) \geq c\, r_ \alpha ^{N-2} G(x_0, A_ \alpha )\, G(A_ \alpha ,x)$ for $x \in   U_ \alpha $. Thus, 
\begin{align}
G({\mathbf 1}_{U_ \alpha }K_0^VV)(x_0)  &= \int_{U_ \alpha } G(x_0, x)\, K_0^V(x)\, V(x)\, dx  \nonumber \\
&\geq { \frac {  1} {c }}\,\,r_ \alpha ^{2N-4}G(x_0,A_ \alpha )\, K_0(A_ \alpha )\,\int _{U_ \alpha }\,    G(A_ \alpha ,x)  \, G^V(A_ \alpha ,x)\, V(x)\,dx \nonumber\\
& \geq { \frac {  1} {c }}\,\,r_ \alpha ^{N-2}\,\int _{U_ \alpha }\,  G(A_ \alpha ,x) \, G^V(A_ \alpha ,x)\, V(x)\,dx
\end{align}
where for the last line,  we have also used Lemma \ref{GKlemma} and Harnack inequalities.

 Similarly (working now only with $V=0$), we have
\begin{equation}K_ 0 (x) \leq c\, r_ \alpha ^{N-2} K_0(A_ \alpha )\,G(A_ \alpha , x), \;x \in   U_ \alpha \end{equation}
and $G(x_0,x) \leq  c\, r_ \alpha ^{N-2} G(x_0, A_ \alpha )\, G(A_ \alpha , x)$ for $x \in   U_ \alpha $. Thus
\begin{align}
G({\mathbf 1}_{V_ \alpha }K_0V)(x_0)  &= \int_{V_ \alpha } G(x_0, x)\, K_0(x)\, V(x)\, dx  \nonumber \\
 & \leq c\,r_ \alpha ^{N-2}\,\int _{V_ \alpha }\,  \vert  G(A_ \alpha ,x) \vert  ^2\, V(x)\,dx.
\end{align}

It follows, on using Proposition \ref{propmain}  (and Harnack inequalities), that 
\begin{align} G({\mathbf 1}_{V_ \alpha }K_0V)(x_0)  & \leq c\, G({\mathbf 1}_{U_ \alpha }K_0^V\,V)(x_0)  
\end{align} 
for every $ \alpha \in   I $ and a constant $c$ independent of $ \alpha $. Summing over $ \alpha $ and taking into account condition (ii) for the family $ \{ B'_ \alpha  \}_{ \alpha  \in   J}$, we get that for $W= \bigcup    V_ \alpha $ 
\begin{align} G({\mathbf 1}_W\, K_0\, V)(x_0) \leq c\, G(K_0^V\, V)(x_0)< \infty. \end{align}

Since $\Omega \setminus W$ is nontangential in $\Omega $ at $0$, and $K_0 \leq c\, K_0^V$ near $0$ on $\Omega \setminus W$ we obtain the desired result: $G(K_0\, V)(x_0)< \infty $. The proof of Theorem \ref{charac} is complete. $\square$

\section {Almost everywhere regularity}

{\bf 6.1.} We first notice that  Corollary \ref{cormain} (see introduction)  easily follows from Theorem \ref{th1}. Clearly, if the integral $\int_\Omega  K_\nu (x)\, V(x)\,  G(x,x_0)\, dx $ is finite then by Fubini's theorem $\int _\Omega K_y(x)\, G(x,x_0)\,V(x)\, dx $ is finite for $\nu$-almost every $y \in   \partial \Omega $. Since $ \nu =f\mu  $ with $f>0$ this means -- by Theorem \ref{th1} -- that $\mu  $-almost every $y$ is in $ { \mathcal R}\mbox{\rm \footnotesize{eg}}_V(\Omega )$.

On the other hand if $\mu  $-a.e.\ $y \in   \partial \Omega $ is $V$-finely regular, then $y \mapsto \int _\Omega \,K_y(x)\, V(x)\, G(x,x_0)\,$ $ dx$ is finite $\mu  $-a.e. So, if we set $A_n:= \{ y \in   A\,;\, n-1 \leq G(K_y\,V)(x_0) < n\, \}$, $B:=\partial \Omega \setminus  \bigcup    _{n \geq 1} A_n$, $ f=  \mathbf 1_B +\sum _{n \geq 1} { \frac {  1} { n^3}}\mathbf 1_{A_n}$ and $\nu := f\, \mu  $  we have $\mu  (B)=0$ and $G(K_\nu V)(x_0)< \infty $.

\vspace{0.7truecm} 
{\bf 6.2.} Recall from the Martin boundary theory \cite{martin} \cite{naim} that a set  $A \subset   \Omega $ is minimally thin at $y \in   \partial \Omega $ for $ L_V:=\Delta  -V$ iff $^\Omega \! \widehat R_{K_y^V}^A$ is a $ L_V$-potential in $\Omega $ (the r\'{e}duite is taken  with respect to the  $L_V$-potential theory). If $A$ is minimally thin at $y$ and closed in $\Omega $, it is well-known that $u:=K_y^V-^\Omega \!\!R_{K_y^V}^A$ is then minimal $L_V$-harmonic in $ \tilde \Omega =\Omega \setminus A$ \cite{naim}. Thus if moreover $ \tilde \Omega $ is a Lipschitz domain and $x_0 \in    \tilde \Omega $, $u$ is proportional to $ k_y^V$, the $L_V$-Martin function in $ \tilde \Omega $ with  pole at $y$ and such that $ k_y^V(x_0)=1$. So $ k_y^V \leq C\, K_y^V$ for some real $C>0$.

 \begin{lemma} \label{lemmaregeffi}Let $ \tilde \Omega $ be a Lipschitz subdomain of $\Omega $, let $y \in   \partial \Omega  \cap \partial  \tilde \Omega $ and assume that  $\Omega \setminus  \tilde \Omega $ is minimally thin at $y$ (in $\Omega $, with respect to $ \Delta  $). Then $y \in   { \mathcal R}\mbox{\rm \footnotesize{eg}}_V(\Omega )$ if and only if $y \in    { \mathcal R}\mbox{\rm \footnotesize{eg}}_V( \tilde \Omega )$. Here $ { \mathcal R}\mbox{\rm \footnotesize{eg}}_V( \tilde \Omega )$ means $ { \mathcal R}\mbox{\rm \footnotesize{eg}}_{V_1}( \tilde \Omega )$ where $V_1=V_{ \vert \tilde \Omega}   $. 
\end{lemma}

We note that at least one implication seems clear with the probabilistic point of view.

{\sl Proof.} Assume as we may that $x_0 \in    \tilde \Omega $ and let $\nu_y $ be a pseudo-normal for $ \tilde \Omega $ at $y \in   \partial \Omega  \cap \partial  \tilde \Omega $.  By Lemma \ref{GKlemma}  it is clear that $ k_y^V \geq C K_y^V$ and $ k_y:= k_y^0\geq C K_y$ near $y$ on $S_y=(y,y+\nu _y]$,  for some $C>0$. Using the facts reminded before the lemma for $V=0$, we have $ k_y\sim K_y$ (where $\sim$ means ``is in between two positive  constants times") on $S_y$ (near $y$).

Now if $ y$ is $V_1$-finely regular for $\tilde \Omega $, $ k_y^V \leq c k_y$  in $ \tilde \Omega $ for some $c>0$. Thus, by the above  $K_y \geq C\, K_y^V$  along $S_y$ near $y$ for another $C>0$ and $y $ is $V$-finely regular as a boundary point for $\Omega $.

Suppose now that $y  \in   { \mathcal R}\mbox{\footnotesize eg}_V(\Omega )$, that is $K_y \geq { \frac {  1} {C }} K_y^V$ in $\Omega $ for some $C \geq 1$. Since $\Omega \setminus  \tilde \Omega $ is (minimally) thin at $y$, there is a (classical) potential $p$ in $\Omega $ such that $p \geq K_y$ in $\Omega \setminus  \tilde \Omega $ and $p$ must be also a $ \Delta  -V$ potential (a nonnegative $ \Delta  -V$ subharmonic minorant is subharmonic and so vanishes).   Thus $\Omega \setminus  \tilde \Omega $ is also minimally thin at $y$ for $ \Delta  -V$ and  along $S_y$ near $y$ we have : $K_y^V\sim K_y \sim k_y$, $K_y^V\sim k_y^V$ and hence also $ k_y^V \leq C'k_y$. So $y \in   { \mathcal R}\mbox{\footnotesize eg}_{V}( \tilde \Omega )$. $\square$

\vspace{5truemm}

{\bf 6.3} {\sl Lebesgue almost everywhere fine regularity.}
\begin{theorem} \label{th3} Let $A$ be a Borel subset of $   \partial \Omega $. The following are equivalent:
\begin{enumerate}
\item   $H_{N-1}( { \mathcal R}\mbox{\rm \footnotesize{eg}}_V(\Omega ) \cap    A)>0$. 
\item  There exists a compact ${ \mathbb K}   \subset   A $ such that $H_{N-1}({ \mathbb K} )>0$ and \begin{equation} \int_\Omega  \omega ^{ \mathbb K} (x)\, G(x_0,x)\, V(x)\, dx < \infty
\end{equation} 
where $\omega^{ \mathbb K} $ is the harmonic measure of ${ \mathbb K} $ in $\Omega $ and $G$ is Green's function in $\Omega $.  
\item There exist a compact set ${ \mathbb K}  \subset   A$ and  a Lipschitz subdomain $ \tilde \Omega $ of $\Omega $ such that $ H_{N-1}( {\mathbb K} )>0$, $\partial \Omega  \cap \partial  \tilde \Omega ={ \mathbb K} $ and    
$\int_{  \tilde  \Omega }  \delta_\Omega  (x)V(x)\, dx< \infty $. 
\end{enumerate} 
\end{theorem}

Recall that by a well-known result  the harmonic measure $\omega _{x_0}:=\omega ^.(x_0)$ of $x_0$ in $\Omega $ and the Hausdorff measure  $H_{N-1}$ are mutually absolutely continuous on $\partial \Omega $ \cite{Dah}.

 {\sl Proof.} $(i)\Rightarrow(ii)$. It suffices to note that by Theorem  \ref{charac} there exist a positive real  $c$ and  a compact  set ${ \mathbb K}  \subset   A$ such that $H_{N-1}({ \mathbb K} )>0$ and  $\int_\Omega  G(x,x_0)\, K_y(x)\, V(x)\, dx  \leq c$ for $y \in   { \mathbb K} $.   Integrating with respect to $y$ against $\omega _{x_0}$ (the harmonic measure  of $x_0$ w.r. to $\Omega $), we get 
 \begin{align}   
\int_\Omega  G(x,x_0)\,\omega ^{ \mathbb K} (x)\,V(x)\, dx =\int _{ \mathbb K} \int_\Omega  G(x,x_0)\, K_y(x)\, V(x)\, dx \, d\omega _{x_0}(y) \leq c\, \omega _{x_0}({ \mathbb K} )< \infty \nonumber
\end{align} 
where we have used the  identity $\omega ^{ \mathbb K} (.)=\int_{{ \mathbb K} } K_y(.)\, d\omega _{x_0}(y)$. 

$(ii)\Rightarrow (iii)$. With the notations and assumptions  of Section {2.2} we may assume that $f$, $r$, $ \rho $ are such that $T(1) \cap    \Omega=U  $,  $x_0\notin  \overline  U$ -- where $U:=U_f(r, \rho )$ -- and that  ${ \mathbb K}  \subset   \partial U \cap    T({\frac  { 1} {8}})$. 

By Dahlberg's results in  \cite{Dah} there is a compact set  ${ \mathbb K} ' \subset   { \mathbb K} $ and a real $c>0$ such that $H_{N-1}({ \mathbb K} ')>0$ and $G_{x_0}(x) \geq c\,  \delta _\Omega (x)$ in every truncated cone $C_y= \{ (x',x_N) \in   { \mathbb R} ^N\,;\,  \vert  x'-y' \vert  <{\frac  { r} { \rho }}\,  (x_N-y_N),\, \, x_N < {\frac  {  \rho } {2}} \, \}$, $y=(y',y_N) \in   { \mathbb K} '$. 
Thus if $ \tilde \Omega = \bigcup    _{y \in   { \mathbb K} '} C_y $ we get by (ii) that $\int _{ \tilde \Omega }  \delta_\Omega  (x)\, V(x)\, dx< \infty $. And $ \tilde \Omega $ is a Lipschitz domain  (in fact ``Lipschitz starlike" with respect to the point $(0,\dots, 0,{\frac  {  \rho } {4}})$). 

$(iii)\Rightarrow(i)$. It is known (cf \cite{HW1}, \cite{HW2}) that for $\omega _{x_0}$-almost every $y \in  { \mathbb K} $ the set $\Omega \setminus  \tilde \Omega $ is minimally thin at $y$. On the other hand, using again Dahlberg's results in \cite{Dah} we may assume -- after replacing $ \tilde \Omega $ by a smaller subdomain in $\Omega $ -- that  $g(x,x_1) \leq c\,  \delta _\Omega (x)$ in $ \tilde \Omega \setminus B(x_1; {\frac  { 1} {2}}  \delta _{ \tilde \Omega }(x_1))$, where $g$ is Green's function for $ \tilde \Omega $ and $x_1$ some fixed point in $ \tilde \Omega $. Thus   $\int _{ \tilde \Omega }  g(x,x_1)\, V(x) \,dx < \infty $. 

If $ \tilde K$ denotes the Martin kernel for $ \tilde \Omega $ with normalization at $x_1$ and if $\mu  $ is the harmonic measure of $x_1$ in $ \tilde \Omega $, as before we have $ \tilde K_\mu  =\mathbf 1$ and $\int _{\partial  \tilde \Omega } \int _{ \tilde \Omega }\, g(x,x_1)  \tilde K_y(x)\, V(x)\, dx \,d\mu  (y)= \int g(x,x_1)\, V(x)\, dx < \infty $. It follows using Theorem \ref{th1} that $\mu  $-almost every (i.e., $H_{N-1}$-almost every) point $y \in    \partial \tilde \Omega $ is finely regular for $ \tilde \Omega $ and $V_{ \vert   \tilde \Omega }$. By lemma \ref{lemmaregeffi}, $(iii)$ follows. $\square$

\begin{corollary} \label{corregae}Let $A$ be a Borel subset of $\partial \Omega $. Then $H_{N-1}$-almost every point $y \in  A$ is finely regular (with respect to $\Omega $ and $V$) if and only if there exists for $H_{N-1}$-almost every point $y \in  A$ a nonempty  open truncated cone $C_y \subset   \Omega $ with vertex at $y$ and such that $\int _{C_y} V(x)\,  \vert  x -y\vert  ^{2-N}\, dx< \infty $.  
\end{corollary}

{\sl Proof.} We have already seen that the condition is necessary (Theorem A.1 in \cite{VYA} or Remark \ref {remarkApp}  above): if $ y    \in   \partial \Omega $ is a $V$-finely regular boundary point, then $\int _C V(x)\,  \vert  x -y\vert  ^{2-N}\, dx< \infty $ for every  open truncated cone $C  $ with vertex at  $y$ which is strictly inner for $\Omega $ (i.e.\  there exists $ \eta  >0$ such that $ \tilde C= \cup _{x \in   C} B(x,   \eta  \,  \vert  x-y \vert  ) \subset   \Omega $).

Suppose now that the condition in the statement holds. It suffices to show that a compact set ${ \mathbb K}  \subset  A$ such that $H_{N-1}({ \mathbb K} )>0$ contains a $V$-finely regular point. Using standard  argument and a rotation  we may also assume that for some positive reals  $\ell,\, K$ and $\alpha$, the cones $C_y=(y+ \Gamma ) \cap  \{ x\, ;\, x_N< \ell \}  $ where $ \Gamma:   = \Gamma  _K:= \{ (x',x_N)\,;\,  \vert  x' \vert  <Kx_N\, \}$, are such that $ \emptyset  \ne C_y \subset   \Omega $ and $\int _ {C_y} V(x)\,  \vert  x -y\vert  ^{2-N}\, dx\leq   \alpha  $ for $y \in   {\mathbb K} $. Diminishing ${\mathbb K} $, we may also assume that the cones $C'_y=(y+ \Gamma  _{K/2}) \cap  \{ x\,;x_N<\ell/2 \}$, $y \in   { \mathbb K} $, have a nonempty intersection.

Denote $L$ the compact subset of all points $y \in   {\mathbb K} $  such that  $H_{N-1}(B(y,t) \cap {\mathbb K} ) \geq  \beta t^{N-1}$ for $0<t \leq {\rm  { diam }}(\Omega )$ where $ \beta >0$ is chosen so small that $H_{N-1}(L)>0$.

Integrating w.r. to $y$ and letting $ \tilde \Omega $ to denote the union of the cones $C'_y$, $y \in   L$, we obtain by Fubini's theorem,
\begin{equation}   \alpha \, H_{N-1}({\mathbb K} )  \geq \int _{\mathbb K} (\int _{C_y} V(x)\,  \vert  x -y\vert  ^{2-N}\, dx)\, dH_{N-1}(y) \geq c\, \int _{ \tilde \Omega }  \delta _\Omega (x)\,V(x)\, dx
\end{equation} 
where $c$ is some positive constant. We have used the fact that if $x \in   C'_y$, then  $x \in   C_z $ for $z \in   {\mathbb K} $ such that $ \vert  z-y \vert   \leq  \varepsilon _0\,  \vert  x-y \vert  $ where $ \varepsilon _0>0$ is independent of $x$ and $y$. Since $ \tilde \Omega $ is a Lipschitz domain it follows from Theorem \ref {th3} that $L \cap { \mathcal R}\mbox{\footnotesize eg}_V(\Omega )\ne  \emptyset  $. $\square$

\begin{corollary} Suppose that $V$ satisfies for some $c \geq 0$ the Harnack property: \begin{align}\sup_{z \in   B(x;{\frac  {  \delta _\Omega (x)} {4}})}  V(z) \leq c\, V(x)\end{align}
for all $x \in   \Omega $. Then under the assumptions of Corollary \ref{corregae}, $H_{N-1}$--almost every $y \in   A \subset   \partial \Omega $ is $V$-finely regular iff for $H_{N-1}$-almost every $y \in   A$ there exists a pseudo-normal $\nu $ at $y$ for $\Omega $ such that  $\int _0^{ \eta  }\, t\, V(y+t\nu )\, dt < \infty $ for  some $ \eta  >0$ such that $ (y,y+  \eta \nu ] \subset   \Omega $.
\end{corollary}

 \section {Generalizations}

 We consider in this section two uniformly elliptic operators in divergence form acting in a bounded region $\Omega $ (at first not assumed to be Lipschitz) ${ \mathcal L}_0= { \mathcal L} -  \gamma $ and ${ \mathcal L}_1= { \mathcal L}- V $ where ${ \mathcal L}=\sum_{1 \leq i,j \leq N}   \partial _j(a_{ij} \partial _i .)$. The $a_{ij}$, $ \gamma $ and $V$ are  Borel measurable and such that : 
 (i) $a_{ij}=a_{ji}$ for $1 \leq i,j \leq N$ (ii) $c_0^{-1}  \vert   \xi   \vert  ^2 \leq   \sum _{i,j} a_{ij}(x) \xi  _i \xi  _j \leq  c_0\,  \vert   \xi   \vert  ^2$ and (iii) $0 \leq  \gamma (x) \leq V(x) \leq a \delta_\Omega  (x)^{-2}$  for  $x \in   \Omega $, $ \xi   \in   { \mathbb R} ^N$ and some positive constants  $a$ and $c_0$.

{\bf 7.A.} We first extend Lemma \ref{mainlemma} and Corollary \ref{maincoro}. The domain $\Omega $ is assumed to satisfy:  $\int _\Omega { \frac { \vert  f(x) \vert  ^2} { \delta (x)^2 }}\, dx \leq C_H\, \int_\Omega   \vert  \nabla f(x) \vert  ^2\, dx$ for all $   f \in   C_0^1(\Omega )$ and a constant $C_H>0$.

\begin{proposition} \label{mainproposition2}Let  $ \zeta_0  \in  \Omega $ and let 
$w(x)=G^0(x, \zeta_0  )$ where $G^0$ is the ${ \mathcal L}_0$-Green's function in $\Omega $. For every   Borel subset $A$ of $\Omega $,
 \begin{align}  \int _A ({\frac  { w(x)} { \delta _{\Omega } (x)}})^2 \, dx  \leq C\; ^{\hspace{1mm}\Omega,{ \mathcal L}_0}\! \widehat R_w^A( \zeta_0  ) 
 \nonumber
\end{align}  where $C$ is a positive constant depending only on  $N$, $C_H  $, $c_0$ and $a$.
\end{proposition}

The r\'{e}duite $^{\hspace{1mm}\Omega,{ \mathcal L}_0}\!R_w^A\,$ (simply denoted $R_w^A $ in what follows)  is here the infimum of all positive ${ \mathcal L}_0$-superharmonic functions $s$ in $\Omega $ such that $s \geq w$ in $A$. It is ${ \mathcal L}_0$-harmonic in $\Omega \setminus \overline A$  and its l.s.c.\ regularization $\widehat R_w^A $  is superharmonic in $\Omega $, in fact the smallest nonnegative superharmonic  function in $\Omega $ which is larger than $ w$ in $A$ except at most on a polar subset.

{\bf Proof.}  We sketch a variant of the proof of Lemma \ref{mainlemma} showing more generally that if $w=G^0(\mu ) $ is the ${ \mathcal L}_0$-potential of $\mu   \in  { \mathcal M}_+(\Omega)$ then   \begin{align} \label{eq71}  \int _A ({\frac  { w(x)} { \delta _{\Omega } (x)}})^2 \, dx  \leq C\;\int _\Omega \,\widehat R_w^A( \zeta  )\, \,d\mu  ( \zeta  ) 
\end{align} 
--- one may also easily adapt the proof of Lemma  \ref{mainlemma} using $H_0^1(\Omega )$ equipped with the scalar product $ \langle .,. \rangle  _H^0$ defined below. Removing from $A$ a polar set we may  assume that $ \widehat R_w^A={}  R_w^A$ and standard approximation arguments  show that we may also  assume that $A$ is open and  relatively compact in $\Omega$ and moreover that $w$ is bounded (using $w_n=w\wedge n$) and the positive measure $\mu  =-{ \mathcal L}_0(w)$ compactly supported in $\Omega $ (using $w'=R_w^U$, with $U$ open such that $\overline A \subset   U \subset    \subset   \Omega $). 

Then $w \in   H_0^1(\Omega  )$  is the potential of $\mu  $ in $\Omega $  with respect to the Dirichlet scalar product $ \langle  \varphi  , \psi  \rangle^0  _H= \sum _{i,j}\int_\Omega   \, a_{ij} \partial _i\varphi\,  \partial _j  \psi \, dx+\int _\Omega \,  \gamma \,  \varphi  \,   \psi \, dx$.  Moreover $ R_w^A \in   H_0^1(\Omega )$ and $R_w^A$ is the ${ \mathcal L}_0$-potential in $\Omega $ of a positive measure $\nu $ supported by the fine closure of $A$ (so $ R_w^A = G^0(\nu )$). Thus
\begin{align} \int _\Omega  {} \,R_w^A( \zeta  )\, \,d\mu  ( \zeta  )  =\int _\Omega  \, \,G^0(\nu )\, d\mu &=\int _\Omega \, G^0(\mu  )\, d\nu =\int _\Omega \, G^0(\nu )\, d\nu  \nonumber \end{align}
using Fubini and that $G^0(\mu  )=G^0(\nu )$ $\nu $-a.e. The result follows since
\begin{align} \int _\Omega  {} \,R_w^A( \zeta  )\, \,d\mu  ( \zeta  ) &= \langle R_w^A,R_w^A \rangle  _H^0 \nonumber
\\ & \geq  C_H\, \int _\Omega  \,({ \frac {  R_w^A(x)} { \delta _\Omega (x) }})^2\, dx \nonumber 
\\& \geq C_H\, \int _A\, ({ \frac {  w(x)} { \delta _\Omega (x) }})^2\, dx . \;\; \square\nonumber \end{align}

\medskip

\begin{corollary}\label{maincorogene2} {\rm Let $\omega $, $\omega '$ be open subsets of ${ \mathbb  R}^N$ such that $\overline \omega ' \subset   \omega $ and let $x_0 \in   \Omega \setminus \overline \omega $. If $\Omega $ is Lipschitz, there  is a $C=C(\omega ,\omega ',\Omega ,x_0,c_0,a)>0$ such that for every positive ${ \mathcal L}_0$-harmonic function $h$ in $\Omega $ vanishing on $\omega  \cap \partial \Omega $ and every (measurable) $A \subset   \omega' $
\begin{align}  
\int _A\,({ \frac {  h(x)} { \delta _\Omega (x) }})^2\, dx  \leq C\,\,h(x_0)\,\,  ^{\Omega,{ \mathcal L}_0 }\!\widehat R_h^A(x_0)
\end{align}}
\end{corollary}

{\bf 7.B. Relative regularity.}

From now on, $\Omega $ is a bounded Lipschitz domain in ${ \mathbb R} ^N$. For $y \in   \partial \Omega $, denote $K_y:=K_y^{ \gamma }$ (resp. $K_y^V$) the ${ \mathcal L}_0$-Martin function (resp.\ the ${ \mathcal L}_1$-Martin function) in $\Omega $ with pole at $y$, normalized at some fixed $x_0 \in   \Omega $. Let  $ \tilde K_y^V$ be the greatest ${ \mathcal L}_1$-subharmonic minorant of $K_y$ (the latter  is positive ${ \mathcal L}_1$-superharmonic). Again, $ \tilde K_y^V=c(y)\, K_y^V$ for some $c(y) \geq 0$. Set $R=V- \gamma $.

\begin{definition}  A point $y \in   \partial \Omega $ is an $R \vert   { \mathcal L}_0 $-regular boundary point (or an $R$-regular boundary point with respect to ${ \mathcal L}_0$) iff $c(y)>0$ which is denoted $y \in   { \mathcal R}\mbox{\footnotesize eg}_{{ \mathcal L}_1 \vert  { \mathcal L}_0}(\Omega )$. If $c(y)=0$, $y$ is a $R \vert   { \mathcal L}_0$-singular boundary point.
\end{definition}

Proposition \ref{equivdef} extends in a straightforward manner. The boundary point $y $ is in $  { \mathcal S}\mbox{\footnotesize {ing}}_{{ \mathcal L}_1 \vert  { \mathcal L}_0}(\Omega )$ iff $K_y^V/K_y$ has an infinite superior limit along a pseudo-normal at $y$, or iff $K_y^V/K_y \to  + \infty $ nontangentially at $y$. Equivalently,  $\lim G^1(x,x_0)/G^{ 0}(x,x_0)=0$ as $x \to  y$, $x \in   \Omega $ -- where  $G^j$ is ${ \mathcal L}_j$-Green's function in $\Omega $.

\vspace{3truemm}

\begin{theorem}\label{charac2} Let $y \in   \partial \Omega $ and $R=V- \gamma $. The following  are equivalent:
\begin{enumerate}
\item $y$ is an $R \vert   { \mathcal L}_0 $-finely regular boundary point of $\Omega $,
\item  $G^{0}(RK_y^V)(x_0)< \infty $,
\item $K_y^V$ admits an ${ \mathcal L}_0$-superharmonic majorant in $\Omega $,
\item $G^{0}(RK_y)(x_0)< \infty $.
\end{enumerate}
\end{theorem}

As in Section 4 one shows that  $(i)\Rightarrow(ii)\Rightarrow(iii)\Rightarrow(i)$ and $(iv) \Rightarrow (i)$. The argument is exactly the same, after replacing $ \Delta  $ by ${ \mathcal L}_0$, $ \Delta  -V$ by ${ \mathcal L}_1$. The proof is omitted.

\vspace{3truemm} 
{\bf 7.C. Necessity of the  ``$\mathbf{G(VK_y)< \infty} $" condition}

To prove the remaining implication $(i)\Rightarrow (iv)$, we adapt the argument in Section 5, a task which turns out to be quite straightforward. Note first 
the following simple extension of Proposition \ref{propmain}. Using same notations as in Proposition \ref{propmain} for $U=U(f,r, \rho )$, $T(t)=B(0,r)\times (- \rho , \rho )$, and $ \zeta  _1=(0,\dots, 0,{ \frac {  t+t'} {2 }} \rho )$, consider now a positive solution $u$ of ${ \mathcal L}_1(u)=0$ in $W:=U_{t'}$, a positive solution $v$ of ${ \mathcal L}_0(v)=0$ in $ W$, with $u=v=0$ on  $\partial _\#U \cap \overline W$, and set (as before) $R=V-V_0$.

\begin{proposition} \label{propmain2}    There is a   constant $C=C(N,{\frac  {  \rho } {r}} ,{\frac  { t'} {t}}  ,a, { \mathcal L}_0)$ $>0$, such that  
\begin{align}\label{keyeq2}  u(  \zeta  _1)\int_{U_t} \,  \vert  v(x) \vert  ^2\, R(x)\,  dx  \leq Cv( \zeta  _1)\, \int _{ W} &v(x)\,R(x)\, u(x)\, dx\nonumber
\\ &= Cv( \zeta  _1) \int _{ W} v(x)  \,{ \mathcal L}_0u(x)\,\, dx  \end{align}
\end{proposition}

{\bf Proof.} a) Again, enlarging $t$ we may assume  $t>t'/2$ and set $  W_1=U_{(t+2t')/3}$, $B_1:=B( \zeta  _1,{\frac  { (t'-t)r} {20}} )$, $   \tilde W_1=U_{(t+2t')/3}\setminus B_1$.  It suffices to prove  (\ref{keyeq2}) with the integration domain $W$  in the right hand side   replaced by $   \tilde W_1=U_{(t+2t')/3}\setminus B_1$.

If $ \tilde V=\mathbf 1_{W\setminus B_1} V$, $ \tilde  \gamma = \mathbf 1_{W\setminus B_1}  \gamma $, and if $ \tilde u$ (resp.\ $ \tilde v$)  is  a positive $ { \mathcal L}  - \tilde V$-harmonic function (resp. a ${ \mathcal L}  - \tilde  \gamma $ harmonic function) in $W=U_{t'}$ vanishing on $\partial _\#U \cap \overline W$, we see as before that $u$ and $  \tilde  u$ (resp.\ $v$ and $  \tilde  v$) normalized at $ \zeta  _1$ are equivalent in size  in $W_1 \cap     \{ x=(x',x_N)\,;\, x_N \leq t \rho   \}$, and by Harnack inequalities this equivalence extends to $W_1$. 
Thus we may replace  $V$ and $ \gamma $ by $ \tilde V$ and $ \tilde  \gamma $, i.e.\ assume that $V$ and $ \gamma $ are supported in $W\setminus B_1$ and  that $v( \zeta   _1)=u( \zeta  _1)=1$.

\vspace{4truemm}

b) Let $w$ denote the (positive) ${ \mathcal L}_0$-harmonic function in $ W_1=U_{ \frac {  t+2t'} {3 }}$ such that $w=u$ in $\partial   W_1$. Exactly as in Section 5, we see that $w \leq  cv$ in $   W_1$, $c^{-1} \leq w( \zeta  _1) \leq c$ and  $c^{-1} w \leq v \leq cw$ in $U_t$, for some  $c=c(N,a,r/ \rho, t'/t ,{ \mathcal L})$.

c) Let   $s(x)=w(x)-u(x)$, $x \in   \overline  {   W_1}$. Now  $s=  G^0(Ru)= G^1(Rw)$ where $R=V-V_0$, $  G^0$ (resp.\ $  G^1$) is Green's function for $    W_1$ with respect to $ { \mathcal L}_0  $ (resp.\ with respect to ${ \mathcal L}_1$). Thus  
\begin{align}
 s( \zeta  _1)  & =  \int _{   W_1 } \,   G^0( \zeta  _1,x)\,\, R(x)\, u(x)\, dx  \leq C\, r^{2-N}\,  \int_{ \tilde W_1} v(x)\, R(x) \,u(x)\, dx\end{align}

since $G^0( \zeta  _1,x) \leq C\, r^{2-N}\, {\frac  { v(x)} {v( \zeta  _1)}}$ in $ \tilde W_1=W_1\setminus B_1$, with $C=C(N,{ \frac {  t'} {t }},{ \frac {   \rho } {r }},a, { \mathcal L})$.

Now  if $A= \{ x \in   U_t\,;\, u(x) \leq {\frac  { 1} {2}} w(x)\,  \}$, $s \geq {\frac  { 1} {2}} w$ in $A$, we have $\, ^{W_1} \! R_w^A \leq 2s$ in $W_1$ (r\'{e}duite with respect to ${ \mathcal L}_0$) and by  Corollary \ref{maincorogene2}
\begin{equation} \int _ A ({\frac  { w(x)} { \delta_W  (x)}})^2 \, dx  \leq 2 \,c\,r^{N-2}\,s( \zeta  _1)  \leq 2\, c'\,\int _{W_1} v(x)\, R(x)\,u(x)\, dx \nonumber
\end{equation}
with $c'=c(N,{ \frac {  r} { \rho  }},{ \frac {  t'} {t }})$. So, since $V \in   { \mathcal V}(a, U)$,  
\begin{align}\int_A w(x)^2\, V(x)\, dx \leq a \int _ A ({\frac  { w(x)} { \delta_W (x)}})^2 \, dx\leq 2\,a\, c'\,\int _{ W_1}\, v(x)\, R(x)\, u(x)\, dx.\end{align} 

On the other hand, since  $w \leq 2u$ on $U_t\setminus A$,
 \begin{align} \int _{U_t\setminus A} (w(x))^2 \, R(x)\, dx&= \int  _{U_t\setminus A} w(x)\, ({\frac  { w(x)} {u(x)}} )\;   { \mathcal L}_0( u)(x)\, dx \nonumber
\\ & \leq  2\int _{W_1} w(x)\, { \mathcal L}_0u(x)\, dx=2\, \int _{W_1}w(x)\,R(x)\, u(x) dx. 
\end{align}  

And (\ref{keyeq2}) follows. $\square$

 To complete the proof of Theorem \ref{charac} we need to show that if $y_0 \in   \partial \Omega $  is a $R \vert  { \mathcal L}_0$-finely  regular boundary point of $\Omega $ -- so  $G^0(RK_{y_0}^V)< \infty $ and $K^V_{y_0} \leq c\, K_{y_0}$ -- , then $G^0(RK_{y_0})(x_0)< \infty $. The argument for the case ${ \mathcal L}= \Delta  $, $ \gamma =0$ treated in Section 5 can be once more repeated. We sketch it for the reader's convenience.

We 	assume  $y_0=0$, and   $\Omega  \cap  \{ (x',x_N) \in   { \mathbb  R}^N\,;\,  \vert  x' \vert  <r,\,  \vert  x_N< \rho \,  \}= U_{f}(r, \rho )$ for a ${ \frac {   \rho } {10 r }}$-Lipschitz function $f:{ \mathbb  R}^{N-1} \to  { \mathbb  R}$ such that $f(0)=0$. 
Using the same   covering $\{B_ \alpha \}_{ \alpha  \in   J} $ of $B'= \{ x' \in   { \mathbb  R}^{N-1}\,;\, 0< \vert  x' \vert   \leq { \frac {  r} {10 }}\, \}$ in ${ \mathbb  R}^{N-1}$ as in the final part of Section 5, and the same related sets  $U_ \alpha =T'_ \alpha  \cap \Omega =T'_ \alpha  \cap U$, $V_ \alpha =T_ \alpha  \cap U$ and points $A_ \alpha$, we get now 
by the (extended)  boundary Harnack principle (Lemma \ref{PHF} and Remark  \ref{remark8}) and  Harnack inequalities, 
\begin{align}
G^0({\mathbf 1}_{U_ \alpha }K_0^VR)(x_0)  &= \int_{U_ \alpha } G^0(x_0, x)\, K_0^V(x)\, R(x)\, dx  \nonumber \\
&\geq { \frac {  1} {c }}\,\,r_ \alpha ^{2N-4}G^0(x_0,A_ \alpha )\, K_0(A_ \alpha )\,\int _{U_ \alpha }\,    G^0(A_ \alpha ,x)  \, G^1(A_ \alpha ,x)\, R(x)\,dx \nonumber\\
& \geq { \frac {  1} {c }}\,\,r_ \alpha ^{N-2}\,\int _{U_ \alpha }\,  G^0(A_ \alpha ,x) \, G^1(A_ \alpha ,x)\, R(x)\,dx
\end{align}
and as $K_ 0 (x) \leq c\, r_ \alpha ^{N-2} K_0(A_ \alpha )\,G^0(A_ \alpha , x),$
and $G^0(x_0,x) \leq  c\, r_ \alpha ^{N-2} G^0(x_0, A_ \alpha )\, G^0(A_ \alpha , x)$ for $x \in   U_ \alpha $, \begin{align}
G^0 ({\mathbf 1}_{V_ \alpha }K_0V)(x_0)  &= \int_{V_ \alpha } G^0(x_0, x)\, K_0(x)\, R(x)\, dx  \nonumber \\
 & \leq c\,r_ \alpha ^{N-2}\,\int _{V_ \alpha }\,  \vert  G^0(A_ \alpha ,x) \vert  ^2\, R(x)\,dx.
\end{align}So that by  Proposition \ref{propmain2}  (and Harnack inequalities) 
\begin{align} G^0({\mathbf 1}_{V_ \alpha }K_0V)(x_0)  & \leq c\, G^0({\mathbf 1}_{U_ \alpha }K_0^V\,V)(x_0)  
\end{align} 
for every $ \alpha $ and a constant $c$ independent of $ \alpha $, and one conclude as in Section 5. 
 $\square$
  
  \medskip
  
  {\bf 7.D. Regularity  $\mathbf{\omega _{x_0}}$-almost everywhere} 
  
  Finally we mention  an extension of the results in Section {\bf 6}. Let $x_0$ be a fixed point in the bounded  Lipschitz domain $\Omega $ and let ${ \mathcal L}_j$ be as above. Assume moreover that $ \gamma \equiv 0$ (so ${ \mathcal L}_0= { \mathcal L}$, ${ \mathcal L}_1={ \mathcal L}_0-V$, $R=V$). Denote $\omega :=\omega _{x_0}$ the harmonic measure of $x_0$ in $\Omega $ with respect to ${ \mathcal L}$ (in general, this measure  is  not absolutely continuous with respect to $H_{N-1}$ on $\partial \Omega $).
  
  Then Theorem \ref{th3} extends to the present situation if one replaces in its statement ${ \mathcal R}\mbox{\footnotesize eg}_V(\Omega )$ by ${ \mathcal R}\mbox{\footnotesize eg}_{{ \mathcal L}_1 \vert   { \mathcal L}_0}(\Omega )$, the measure $H_{N-1}$ in $\partial \Omega $ by $\omega $ (three times), $G$ by $G^0$ and $ \delta _\Omega (x)$ in $(iii)$ by $g^0(x,x_0)$ -$g^0$ being the ${ \mathcal L}_0$-Green's function in $ \tilde \Omega $-.  The function $\omega _{\mathbb K} $ is now  the ${ \mathcal L}_0$-harmonic measure of ${\mathbb K} $ in $\Omega $. 
  Corollary \ref{corregae} extends then as follows.

\begin{corollary} \label{corregae2} If $A$ is a Borel subset of $\partial \Omega $,  $\omega $-almost every point $y \in  A$ is in ${ \mathcal R}\mbox{\rm \footnotesize eg}_{{ \mathcal L}_1 \vert  { \mathcal L}_0}(\Omega )$ iff for $\omega $-almost every $y \in  A$  there is a nonempty  open truncated cone $C_y \subset   \Omega $ with vertex at $y$ such that $\int _{C_y} V(x)\,  \vert  x -y\vert  ^{2-N}\, dx< \infty $.  
\end{corollary}

The proof for ${ \mathcal L}_0= \Delta  $ in Section 6 is easily adapted provided we set now $L:= \{ y \in   {\mathbb K} \,;\, \omega (B(y,t) \cap {\mathbb K} ) \geq  \beta \omega (B(y,t)$ for $0<t<\text{diam}(\Omega )\, \}$.

\vfill

\pagebreak

\noindent

\end{document}